\documentclass[11pt, reqno]{amsart}
\usepackage{url}
\usepackage[colorlinks]{hyperref}
\usepackage{amssymb,amsmath,amscd,amsthm}
\usepackage{graphicx}
\usepackage{mathtools}
\usepackage[margin=1in]{geometry}
\usepackage{float}
\usepackage{subcaption}

\hypersetup{
     colorlinks=true,
     linkcolor=blue,
     filecolor=blue,
     citecolor = blue,      
     urlcolor=blue,
     }
\usepackage[
    natbib=true,
    style=numeric-comp,
    sorting=none,
    backend = bibtex
]{biblatex}
\usepackage{amsfonts}
\usepackage{bold-extra}
\usepackage[title]{appendix}
\graphicspath{ {figures/} }

\makeatletter
\renewcommand\subsection{\@startsection{subsection}{2}%
  \z@{-.5\linespacing\@plus-.7\linespacing}{.5\linespacing}%
  {\normalfont\bfseries\scshape}}
\renewcommand\subsubsection{\@startsection{subsubsection}{3}%
  \z@{.5\linespacing\@plus.7\linespacing}{-.5em}%
  {\normalfont\bfseries\itshape}}
\makeatother

\newtheorem{theorem}{Theorem}
\newtheorem{definition}{Definition}

\numberwithin{equation}{section}

\newcommand{\pd}{\text{\sc pd}}
\newcommand{\dowker}{\mathcal{D}}
\newcommand{\VR}{\mathcal{VR}}

\DeclareUnicodeCharacter{02B9}{'}
\DeclareUnicodeCharacter{0308}{''}

\begin{document}

\title[Relational persistent homology for multispecies data]{Relational persistent homology for multispecies data \\ with application to the tumor microenvironment}

\author[Stolz, Dhesi, Bull, Harrington, Byrne, Yoon]{Bernadette J. Stolz\textsuperscript{1,2}, Jagdeep Dhesi\textsuperscript{2}, Joshua A. Bull\textsuperscript{2},\\ Heather A. Harrington\textsuperscript{2,3*}, Helen M. Byrne\textsuperscript{2,4*}, Iris H.R. Yoon\textsuperscript{2,5*}}

\maketitle
{\tiny
\noindent
\textsuperscript{1}Laboratory for Topology and Neuroscience, EPFL, Station 8, Lausanne, Switzerland. \\
\textsuperscript{2}Mathematical Institute, University of Oxford, Andrew Wiles Building, Woodstock Rd, Oxford, United Kingdom. \\
\textsuperscript{3}Wellcome Centre for Human Genetics, University of Oxford, Roosevelt Dr, Headington, Oxford, United Kingdom. \\
\textsuperscript{4}Ludwig Institute for Cancer Research, University of Oxford, Old Road Campus Research Build, Roosevelt Dr, Headington, Oxford, United Kingdom. \\
\textsuperscript{5} Department of Mathematics and Computer Science, Wesleyan University, 265 Church Street, Middletown, United States of America.
}

{\tiny \noindent \textsuperscript{*}Corresponding e-mails: harrington@maths.ox.ac.uk, helen.byrne@maths.ox.ac.uk, hyoon@wesleyan.edu}


\begin{abstract}Topological data analysis (TDA) is an active field of mathematics for quantifying shape in complex data. Standard methods in TDA such as persistent homology (PH) are typically focused on the analysis of data consisting of a single entity (e.g., cells or molecular species). However, state-of-the-art data collection techniques now generate exquisitely detailed multispecies data, prompting a need for methods that can examine and quantify the relations among them. Such heterogeneous data types arise in many contexts, ranging from biomedical imaging, geospatial analysis, to species ecology. Here, we propose two methods for encoding spatial relations among different data types that are based on Dowker complexes and Witness complexes. We apply the methods to synthetic multispecies data of a tumor microenvironment and analyze topological features that capture relations between different cell types, e.g., blood vessels, macrophages, tumor cells, and necrotic cells. We demonstrate that relational topological features can extract biological insight, including the dominant immune cell phenotype (an important predictor of patient prognosis) and the parameter regimes of a data-generating model. The methods provide a quantitative perspective on the relational analysis of multispecies spatial data, overcome the limits of traditional PH, and are readily computable.
\end{abstract}

\section{Introduction}\label{Introduction}

Topological data analysis (TDA) is a field of mathematics that develops topological tools for detecting the shape of data. A prominent tool in TDA, persistent homology (PH) \cite{Ghrist_barcodes, PH_survey, Carlsson2009TopologyAD, ELZ_persistence}, constructs a nested sequence of topological scaffolds of shapes from data, called a filtration of simplicial complexes. 
PH examines the evolution of topological features such as connected components (dimension 0) and loops (dimension 1) across the filtration. The filtration is constructed from meaningful aspects of the data at multiple scales such as distances \cite{Vietoris1927berDH, edelsbrunner1993union}, function values \cite{Chazal_persistence_clustering, Gnther2012EfficientCO}, and densities \cite{CarlssonMulti, multiparamIntro, vipond2021multiparameter}. One possible input to PH is point cloud data, and the output is a persistence diagram, which can be vectorized and integrated with statistics and machine learning methods \cite{ali2022survey}. PH provides an automatic, robust, and interpretable method for analyzing data arising in many fields of biology and medicine, including cancer biology \cite{PH_prostatecancer, TDA_histology, PH_glandular_architecture, PH_breastcancer, Nicolau2011TopologyBD, PH_epithelial, Nardini2021,Stolz2020, yang2023topological}, neuroscience \cite{GridCellTorus, PlaceCellsCurtoItskov, DabaghianSpatialMap, CliqueTop}, and genomics \cite{PH_cancer_gene_network, Cmara2017TopologicalMF, homologous_proteins, emmett2015multiscale, top_evolution}.

Most existing PH applications are limited to the study of data relating to one species. Advanced data collection techniques now generate multispecies data in which distinct species may interact. Data of this nature are ubiquitous in science, ranging from cancer biology and ecology to geospatial analysis. By studying the spatial relationships among species, we can glean insights that would otherwise be missed in non-spatial analyses. 
Extracting spatial relationship information from such data, therefore, requires the development of novel analysis techniques. 
Recently, two topological methods have been proposed to study multispecies data \cite{Dhananjay_heterogeneous, chromaticAlpha}. The first approach concatenates topological features from different cell types in cancer images \cite{Dhananjay_heterogeneous} but does not capture spatial relations between the different cell types. Another method, the chromatic Alpha complex \cite{chromaticAlpha}, encompasses relations among species by constructing a multispecies version of the Delaunay triangulation; its computational implementation and interpretation are still under development.

Here, we present two topological approaches for encoding spatial relations among different species directly at the input level for PH. We implement and showcase these methods on synthetic multispecies data generated by an agent-based model (ABM) of the tumor microenvironment. We show that topological relations encode biological insight by predicting the dominant immune cell phenotype and by clustering the parameter regimes of the data-generating model using the relational topological features. 

Mathematically, the multispecies data we consider can be viewed as a labeled point cloud $P = \bigcup_{i=0}^{m} P_i$ that consists of $m+1$ different species whose spatial distributions may be related to one another. Each point $p \in P$ is in $\mathbb{R}^2$ \footnote{Both methods can be applied to point clouds in $\mathbb{R}^n$ for $n \geq 2$.}. We generated synthetic multispecies spatial data from an ABM that simulates the behavior of different cell types in a tumor microenvironment \cite{Bull2023}. The proposed topological methods are built on Dowker complexes \cite{dowker1952} and witness complexes \cite{deSilva2004}. These \emph{relational PH} methods, which we refer to as \emph{Dowker PH} and \emph{multispecies witness PH}, use one species, e.g., $P_0$, as the potential vertex set for a simplicial complex and use another species to create a filtration.

Dowker PH \cite{Chowdhury2018} is based on a Dowker complex \cite{dowker1952}, which is a simplicial complex that represents relations between two point clouds. Dowker complexes have been used to capture relations in molecular biology \cite{dowker_protein}, networks \cite{Chowdhury2018}, PDF parsers \cite{dowker_metric_comp}, and persistence diagrams \cite{analogous_bars}. We propose using Dowker PH \cite{Chowdhury2018}, a natural extension of Dowker complexes, for multispecies data. Dowker PH of the pair $(P_i, P_j)$ creates a filtered Dowker complex on points $P_i$ based on proximity to points in $P_j$. Dowker PH then examines the topological features of the Dowker complex that evolve as one varies the distances between $P_i$ and $P_j$. The resulting Dowker persistence diagram is agnostic to the choice of $P_i$ or $P_j$ as the vertex set and can informally be interpreted as capturing shared topological features, i.e., connected components and loops, between $P_i$ and $P_j$.  

While Dowker PH encodes pairwise relations, it does not capture how one species, say $P_0$, relates to all other species in $P$. To capture differences between all relations among every pair $(P_0, P_i)$, we present a second approach called \emph{multispecies witness PH}, which is inspired by the lazy witness filtration \cite{deSilva2004}. The multispecies witness filtration first creates a Delaunay triangulation \cite{Delauney1934} on $P_0$ and creates a filtration based on the number of points in $P_i$ close to simplices in $P_0$. We chose the Delaunay triangulation because of its simplicity and close relationship to the lazy witness filtration (see Theorem 3 in \cite{deSilva2004}) \footnote{Note that our construction differs from the lazy witness filtration where the filtration values of the simplices in $P_0$ are determined by their proximity to witnesses.}. To encode $P_0$'s relation to all other subpopulations, we construct $m$ separate filtrations, measure the distance between their topological features and combine these distances into a \emph{topological distance vector}, which can then be used as input into classification or machine learning tools.

The paper is organized as follows. In Section \ref{sec:data}, we describe the synthetic multispecies data and introduce the two  questions arising in the study of data from the tumor microenvironment. In Section \ref{sec:single system TDA}, we briefly review the mathematical preliminaries of PH. In Section \ref{sec:multisystem PH}, we present the relational PH approaches designed for capturing relations among multiple species: Dowker PH and multispecies witness PH. In Section \ref{sec:results}, we showcase these methods on a simulated tumor microenvironment and address the biologically motivated questions introduced in Section \ref{sec:data}. The paper concludes in Section \ref{sec:discussion} where we discuss our results and outline directions for future research.

\section{Multispecies spatial data}
\label{sec:data}
We introduce the data set we later analyze, which is synthetic point clouds of multiple species in a simulated tumor microenvironment. Next, we state the two associated domain-specific questions that motivate this mathematical study.

\subsection{Point clouds simulated via agent-based modeling}
\label{sec:ABM}
We study point clouds representing a dynamic and spatially-resolved tumor microenvironment generated by an agent-based model (ABM). ABMs simulate the emergent behavior of a system through the enactment of rules that determine the outcome of interactions between their constituent `agents', here typically individual cells \cite{ABM_Bonabeau}. They are ideally suited to create multispecies data. We use the ABM presented in \cite{Bull2023}. See Appendix \ref{appendix:ABM} and \cite{Bull2023} for details.

Each simulation produces a point cloud $P$ consisting of five species $P = P_T \cup P_S \cup P_N \cup P_M \cup P_V$. Each labeled point cloud represents the locations of tumor cells ($P_T$), stromal cells ($P_S$), necrotic cells ($P_N$), macrophages ($P_M$), and blood vessels ($P_V$). 
The spatial locations of the blood vessels are randomized at the start of each simulation and then held fixed. By contrast, all other cell types are assumed to be motile. Their movement is determined by interactions among the cells and five different diffusible species (oxygen, CSF-1, TGF-$\beta$, CXCL12, and EGF). We focus on simulations that arise by varying two key parameters of the model that affect the behavior of macrophages:
$\chi^m_c$, the chemotactic sensitivity of macrophages to spatial gradients of one of the chemical species (CSF-1), and $c_{1/2}$, a parameter regulating the rate at which macrophage extravasate from the blood vessels~\cite{Bull2023}. We consider 9 different values for each parameter. For each of the 81 possible parameter pairs $(\chi^m_c, c_{1/2})$, we generate up to 20 realizations of the ABM in which the positions of the blood vessels are varied\footnote{These come from 2 sets of 10 realizations in which the threshold value of TGF-$\beta$ required to change macrophage phenotype was varied (either 0.05 or 0.5). Varying this parameter had no qualitative effect on the simulations, and hence the parameter regimes have here been combined.}. Each simulation runs for 500 hours. We focus on the behaviors of macrophages and tumor cells. Each macrophage has an associated phenotype, $\Omega \in [0,1]$, which determines how it interacts with tumor cells. Macrophages with low $\Omega$ have high tumor-killing capacity. Those with high $\Omega$ assist the migration of tumor cells towards the vasculature, thereby promoting metastasis. We refer to macrophages with phenotype $0 \leq \Omega < 0.5$ as $M_1$ or anti-tumor macrophages; we refer to those with phenotype $0.5 \leq \Omega \leq 1$ as $M_2$ or pro-tumor macrophages. 

Simulations are initially seeded with a small cluster of tumor cells at the center of the domain, with blood vessels clustered around the edge. Blood vessels act as sources of oxygen, which is consumed by both stromal cells and tumor cells. Tumor cells are sources of CSF-1, which diffuses through the domain and acts as a stimulus for the recruitment of macrophages and as a chemoattractant for them. 
During each simulation, macrophages with phenotype $\Omega=0$ enter the domain at a rate determined by CSF-1 levels at the blood vessels, with higher CSF-1 increasing the rate of macrophage extravasation. As a macrophage migrates through the domain, its phenotype changes in response to local levels of the different chemical species, including TGF-$\beta$. (For details, see Appendix Section \ref{appendix:ABM}). 

For a given parameter set, at the end of each simulation ($t=500$ hours), we observe one of three distinct qualitative behaviors: 
\begin{itemize}
    \item \textbf{tumor elimination}, in which $M_1$ macrophages dominate the simulation and the tumor cells have been eliminated;
    
    \item \textbf{tumor equilibrium}, in which macrophages are unable to eliminate the tumor cells which form a compact mass, surrounded by macrophages that are predominantly of an $M_1$ phenotype;
    
    \item \textbf{tumor escape}, in which $M_2$ macrophages enhance tumor cell migration to the vasculature. These simulations are characterized by the formation of perivascular niches in which $M_2$ macrophages, tumor cells, and blood vessels are found in close proximity. Such behavior is associated with metastasis of tumor cells \cite{Arwert2018}.
\end{itemize}

We consider two subsets of data generated by the ABM. The first data subset is generated from $2$ realizations of $9\times 9$ parameter combinations of $c_{1/2}$ and $\mathcal{X}_c^m$. The point clouds are generated at $6$ time points ($t = 250, 300, 350, 400, 450, 500$ hours) of the simulation, resulting in $972 = 6 \times 2 \times 9 \times 9$ point clouds. For the second data subset, we consider up to $20$ realizations of $9\times 9$ parameter combinations of $c_{1/2}$ and $\mathcal{X}_c^m$, i.e., a maximum of $1620$ point clouds. As noted in \cite{Bull2023}, limitations on HPC time meant that for some parameter combinations, fewer than 20 realizations were available, giving a total of 1485 point clouds generated at a single `endpoint' time ($t = 500$ hours). For each point cloud, we use the positions of tumor cells, blood vessels, and macrophages (with and without knowledge of macrophage phenotype) as input. For comparison, we also construct simple, i.e., non-topological, descriptor vectors with entries
corresponding to the number of tumor cells, the number of macrophages, the number of
necrotic cells, the average distance of tumor cells to the nearest blood vessel, the average
distance of necrotic cells to the nearest blood vessel, and the average distance of macrophages
to the nearest blood vessel.

\subsection{Statement of biologically motivated problems}
We address the following two biologically motivated questions regarding macrophage and tumor behavior: 
\begin{enumerate}
    \item Can relational PH predict the dominant macrophage phenotype from the cell locations without knowledge of the phenotypes of individual macrophages?
    \item Can relational PH identify the parameter regimes of the ABM that lead to different qualitative behaviors: tumor elimination, escape, and equilibrium with macrophages? 
\end{enumerate}

\begin{figure}[h!]
\centering
\includegraphics[width=0.9\textwidth]{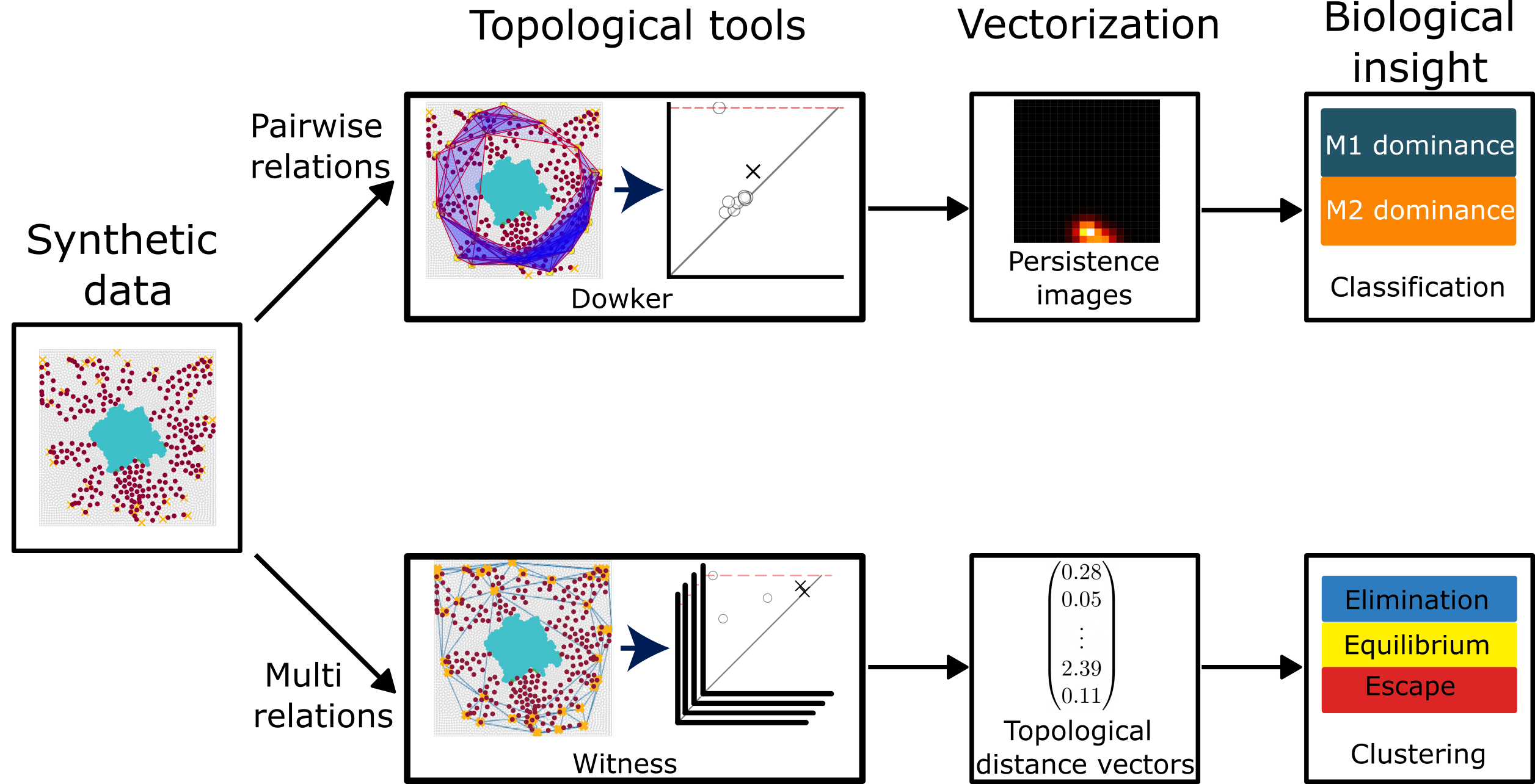}
\captionsetup{width=1\textwidth}
\caption{\textbf{Pipeline and analysis.} We use point clouds generated by an ABM as input to two different topological methods for encoding relations: Dowker PH and multispecies witness PH. We vectorize Dowker topological descriptors using persistence images and vectorize witness topological descriptors via distances between them. Finally, we perform supervised binary classification to predict the dominant macrophage phenotype using Dowker features and perform unsupervised clustering to infer the parameter regimes of elimination, equilibrium, and escape using multispecies witness features.}
\label{fig:full pipeline}
\end{figure}

These two questions motivated the two different pipelines shown in Fig.~\ref{fig:full pipeline}, with the first problem corresponding to the proposed pipeline in the top row and the second problem corresponding to the pipeline introduced in the bottom row.

\subsubsection*{Problem 1: prediction of dominant macrophage phenotype}

\begin{figure}[H]
    \centering
    \captionsetup{width=1\textwidth}
        \includegraphics[width=0.8\textwidth]{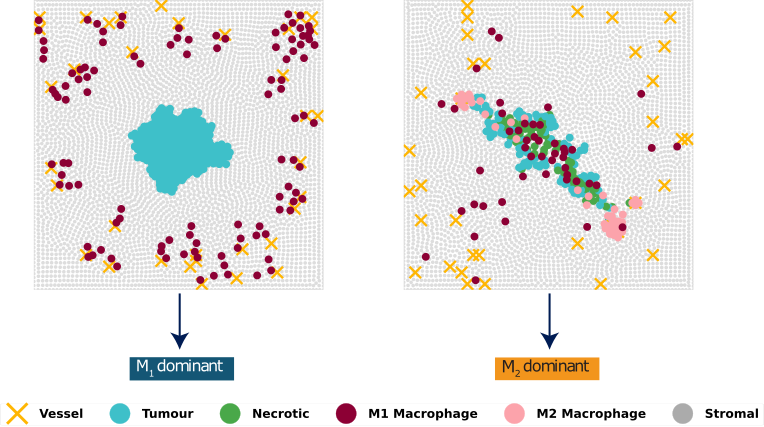}
    \caption{{\bf Problem 1: Prediction of dominant macrophage phenotype.} Given a simulated tumor microenvironment, can we predict the dominant macrophage phenotype?
    }
    \label{fig:joshlabelsI}
\end{figure}

We examine whether relational features can predict the dominance of $M_1$ and $M_2$ macrophages (see Fig. \ref{fig:joshlabelsI}), which is an important predictor of a cancer patient's overall survival time \cite{macrophage_phenotype_survival}. Macrophage phenotype prediction problems may arise in experimental and clinical settings when analyzing imaging data that contains a single macrophage marker or when conventional time- and resource-intensive methods of characterizing macrophage phenotype are not viable \cite{Imaged_ML_macrophage, immunohistochemistry1, immunohistochemistry2, flowcytometry}. We use Dowker PH for this task due to its pairwise encoding of relations. Dowker's shared topological features allow biological interpretation of which relative cell locations directly influence macrophage phenotype. We demonstrate that relational PH can identify the dominant macrophage phenotype based on the spatial relations among the constituents.

\subsubsection*{Problem 2: classification of parameter regimes leading to different qualitative behaviors of the ABM}

\begin{figure}[H]
    \centering
    \captionsetup{width=1\textwidth}
        \includegraphics[width=\textwidth]{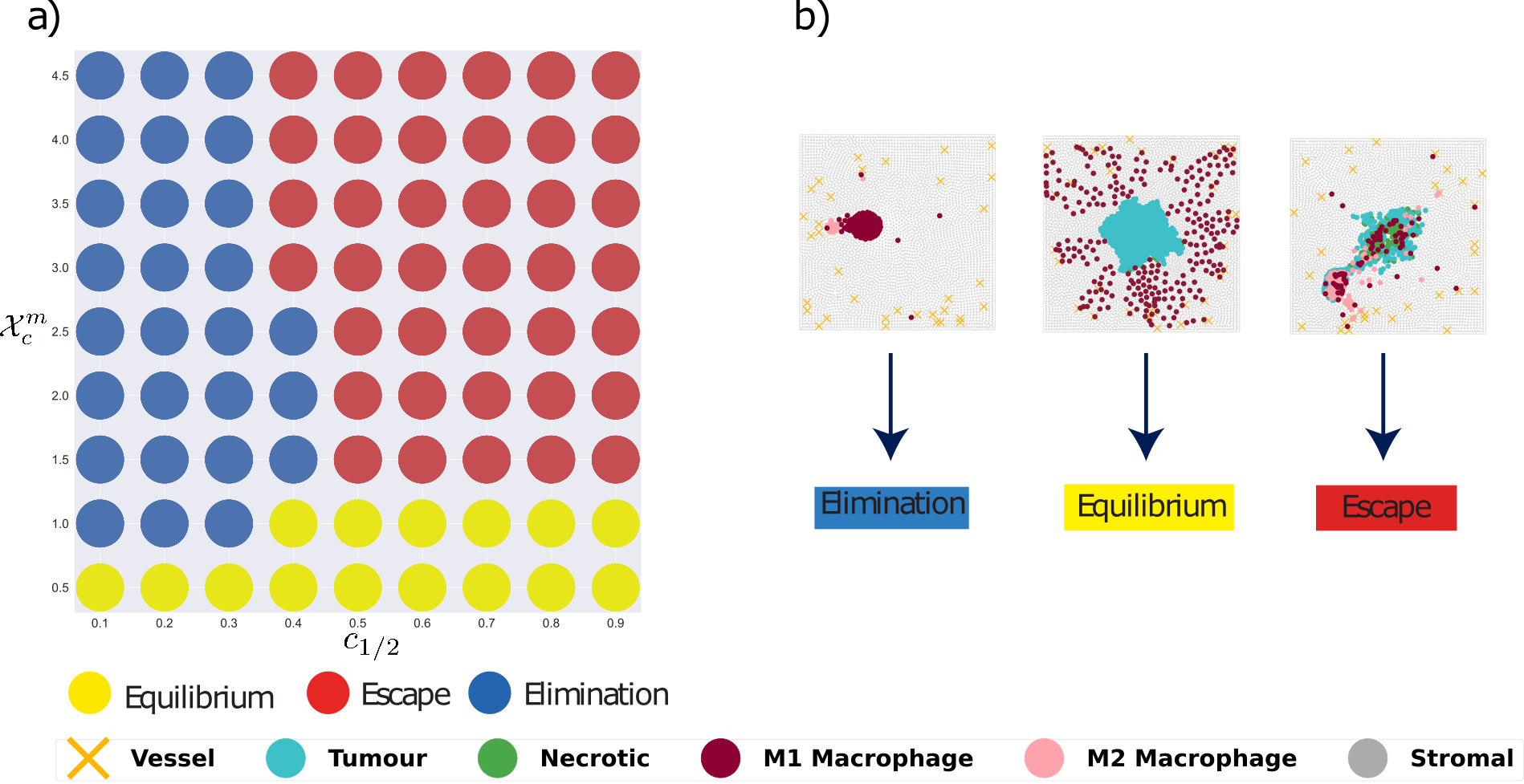}
    \caption{{\bf Problem 2: Classification of parameter regimes leading to different qualitative behaviors of the ABM.} a) Parameter values of $c_{1/2}$ and $\mathcal{X}_c^m$ varied in the ABM. Depending on the parameter combination, a simulation of the tumor microenvironment results in one of three qualitative behaviors: elimination of the tumor (blue), equilibrium of tumor cells and macrophages (yellow), and escape of the tumor cells towards blood vessels (red). The parameter combinations are colored according to the subjective classification of the qualitative behavior observed in one simulation of the model. b) Can we systematically determine the different qualitative behaviors of the ABM from the locations of the different cell types? 
    }
    \label{fig:joshlabelsII}
\end{figure}

Secondly, we explore the use of relational PH in understanding the parameter regimes used to generate different simulations, specifically to classify different parameter regimes from the spatial distribution of the different cell types (see Fig. \ref{fig:joshlabelsII}). The ABM parameters influence the spatial distributions of different cell types in the tumor microenvironment, leading to different tumor compositions and morphology. The qualitative behaviors\footnote{The qualitative behaviors were subjectively assigned in \cite{Bull2023}.} that arise from the different parameter combinations of the ABM are shown in Fig. \ref{fig:joshlabelsII} a). Capturing these differences objectively from the spatial patterns of cells could pave the way for the automated identification of disease stages in microscopy images. Since we are interested in classifying long-term tumor outcomes (escape, elimination, and equilibrium), we consider the ABM output at a single late `endpoint' time ($t = 500$ hours) for varying combinations of parameters $c_{1/2}$ and $\mathcal{X}_c^m$. Multispecies witness PH is ideally suited to this task since it simultaneously takes into account all species in the data set and focuses on their differences.

\section{Mathematical Preliminaries}
\label{sec:single system TDA}

We briefly introduce the standard PH, which can be used to analyze the spatial patterns of point cloud data. For details of PH, see \cite{Ghrist_barcodes, PH_survey, Carlsson2009TopologyAD, ELZ_persistence}.

\subsection{Persistent homology}
Let $P$ denote a point cloud of data in $\mathbb{R}^n$. Here, $P$ is a point cloud of data in $\mathbb{R}^2$ describing the spatial location of biological cells such as cancer cells. The spatial patterns and structure of $P$ can be studied by constructing filtered simplicial complexes, i.e., collections of vertices, edges, triangles, and their higher-order counterparts that can be glued together to approximate topological spaces. We refer to each building block as a simplex. A 0-simplex is a single point in $P$, a $1$-simplex is an edge between two points in $P$, a $2$-simplex is a triangle among three points, and so on. We denote an $n$-simplex by the collection of $n+1$ vertices $(p_0, \dots, p_n)$ that are involved. The standard choice of a filtered simplicial complex is the Vietoris-Rips filtration $\VR_P$ \cite{Vietoris1927berDH}:

\begin{definition}[Vietoris-Rips filtration]
\label{def:VR complex}
Let $P$ be a point cloud and let $d$ be a distance function among $P$. The \emph{Vietoris-Rips complex at parameter $\varepsilon$}, denoted $\VR_P^{\varepsilon}$, is a simplicial complex that has $P$ as the vertex set and has the $n$-simplex $\sigma = (p_0, \dots, p_n)$ if $d(p_i, p_j) \leq \varepsilon$ for all $p_i, p_j \in \sigma$. A \emph{Vietoris-Rips filtration} $\VR_P^{\bullet}$ is a nested sequence of simplicial complexes $\VR_P^{\varepsilon}$ for varying $\varepsilon$.
\end{definition}

\begin{figure}
    \centering
    \captionsetup{width=1\textwidth}
    \includegraphics[width=\textwidth]{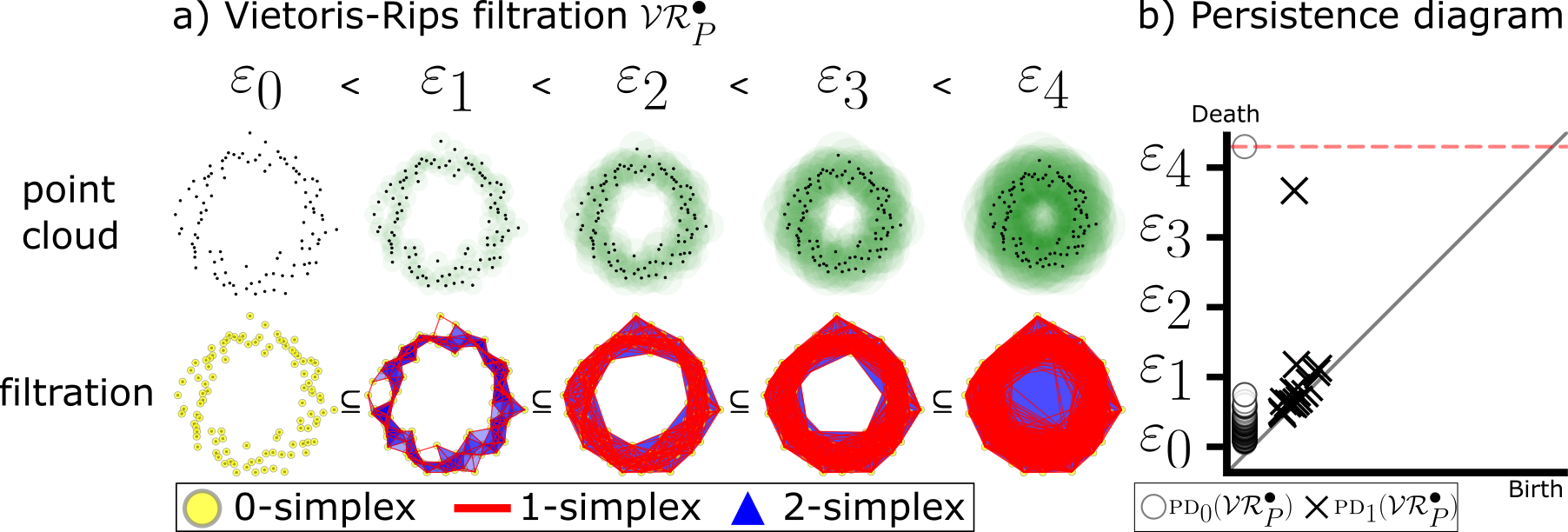}
    \caption{\textbf{An example Vietoris-Rips filtration.} a) Example Vietoris-Rips complexes $\VR_P^{\varepsilon}$ at various $\varepsilon$ parameters. The top row shows the point cloud (in black) and $\varepsilon/2$-ball neighborhoods around each point (in green) for varying $\varepsilon$ values. The bottom row shows the Vietoris-Rips filtration. For a fixed $\varepsilon$, 
   we draw a 0-simplex for each point in the point cloud. Whenever the green balls intersect, we place a 1-simplex between the two corresponding points. We then fill in any higher-dimensional simplices that arise. b) A persistence diagram provides a visual summary of the evolution of connected components (dimension 0, denoted $\pd_0$) and loops (dimension 1, denoted $\pd_1$). We show an overlay of the dimension-0 persistence diagram $\pd_0(\VR^\bullet_P)$ (in circle) and dimension-1 persistence diagram $\pd_1(\VR^\bullet_P)$(in cross). A point on the persistence diagram represents a topological feature. The $x$-coordinate is the parameter $\varepsilon$ at which the feature is born, and the $y$-coordinate is the parameter at which the feature dies. In $\pd_0$, all connected components share the same birth parameter, and the death of a component occurs when two components merge. The red line indicates an infinite death value. There is one point with an infinite death parameter, indicating that the Vietoris-Rips filtration has a single connected component that never vanishes as we increase $\varepsilon$. In $\pd_1$, there is a single point far from the diagonal, indicating that there is one significant loop with a small birth parameter and large death parameter. The remaining points can be considered as noise.}
    \label{fig:vr_filtration}
\end{figure}

The Vietoris-Rips complex $\VR_P^\varepsilon$ at parameter $\varepsilon$ represents the connectivity of $P$ up to proximity $\varepsilon$ (see Fig. \ref{fig:vr_filtration}a). The Vietoris-Rips filtration $\VR_P^{\bullet}$ encodes the connectivity of the point cloud at various proximity parameters. PH provides the means to study topological features such as connected components ($H_0$) and cycles ($H_1$) across nested simplicial complexes. Throughout this paper, we fix the field $\mathbb{F} = \mathbb{Z}/2 \mathbb{Z}$. 

\begin{definition}[Persistent homology]
\label{def:PH} 
Given a nested sequence of simplicial complexes
\[ X^\bullet =X^1 \xhookrightarrow{\iota^1} X^2 \xhookrightarrow{\iota^2} \cdots \xhookrightarrow{\iota^{N-2}} X^{N-1} \xhookrightarrow{\iota^{N-1}} X^N, \]
the \emph{dimension-k persistent homology} of $X^\bullet$ is a collection of $\mathbb{F}$-vector spaces 
\[PH_k(X^\bullet) =  H_k(X^1; \mathbb{F}) \xrightarrow{\phi^1} H_k(X^2; \mathbb{F}) \xrightarrow{\phi^2} \cdots \xrightarrow{\phi^{N-1}}  H_k(X^N; \mathbb{F}),\]
with $\phi^{\varepsilon}$ being the maps induced by $\iota^{\varepsilon}$. 
\end{definition}

The evolution of structural features across a filtration is obtained via the structure theorem.

\begin{theorem}[\cite{carlsson2005} Structure Theorem for persistent homology]
\label{theorem:structure theorem}
Any dimension-$k$ persistent homology $PH_k(X^\bullet)$ obtained from a finite filtered simplicial complex $X^\bullet$ decomposes uniquely as
\[ PH_k(X^\bullet) \cong \bigoplus_{i} I_{b_i, d_i},\]
where each $I_{b_i, d_i}$, called \emph{an interval module}, is a sequence of $\mathbb{F}$-vector spaces
\[I_{b_i, d_i} = 0 \xrightarrow{\phi^0} \cdots \xrightarrow{\phi^{b_i-1}} \mathbb{F} \xrightarrow{\phi^{b_i}} \cdots \xrightarrow{\phi^{d_i-1}} \mathbb{F} \xrightarrow{\phi^{d_i}} 0 \xrightarrow{\phi^{d_i + 1}} \cdots \xrightarrow{\phi^{N-1}} 0\]
with $\phi^{\varepsilon}$ as identity maps for $\varepsilon \in [b_i, d_i)$ and zero otherwise. 
\end{theorem}

Given an interval module $I_{b_i, d_i}$, the parameters $b_i$ and $d_i$ are referred to as the \emph{birth} and \emph{death} times of $I_{b_i, d_i}$. The length (death - birth) is referred to as \emph{persistence}. The decomposition of $PH_k(X^\bullet)$ is often represented using the collection of birth and death times, and they are visualized using a persistence diagram (see Fig. \ref{fig:vr_filtration}b). We denote the dimension-$k$ persistence diagram by $\pd_k(X^\bullet).$ 

Persistence diagrams are stable \cite{stability_Chazal}. That is, there exist distances on persistence diagrams such that small perturbations of the input $P$ result in small changes in the persistence diagram. Two commonly used distances on persistence diagrams are the Wasserstein distance \cite{Wasserstein} and the bottleneck distance \cite{bottleneck}, which are described as follows.

\begin{definition} 
Given two points $x = (x_b, x_d)$ and $y = (y_b, y_d)$ in a persistence diagram let $\| x - y \|_{\infty} = \max \{ \lvert y_b - x_b \rvert, \lvert y_d - x_d \rvert \}$.  Given two persistence diagrams $\pd_k(X^{\bullet})$ and $\pd_k(Y^\bullet)$, the \emph{$q-$Wasserstein distance} is 
\[d_W(\pd_k(X^\bullet), \pd_k(Y^\bullet)) = \inf_{\gamma: \pd_k(X^\bullet) \to \pd_k(Y^\bullet)} \left( \sum_{x \in \pd_k(X^\bullet)} \| x - \gamma(x) \|^q_{2} \right)^{1/q},\]
where $\| \cdot \|_{2}$ is the $L_2$ norm\footnote{Note that the Wasserstein distance can be defined for any type of norm, a typical choice is $L_\mathtt{p}$ for $\mathtt{p} \in [1,\infty]$.}, and the \emph{bottleneck distance} is
\[ d_B(\pd_k(X^\bullet), \pd_k(Y^\bullet)) = \inf_{\gamma: \pd_k(X^\bullet) \to \pd_k(Y^\bullet)} \sup_{x \in \pd_k(X^\bullet)} \| x - \gamma(x) \|_{\infty}  .\]
where $\gamma$ denotes a bijection between $\pd_k(X^{\bullet})$ and $\pd_k(Y^\bullet)$.
\end{definition}

In Section \ref{sec:WitnessFeatures}, we use both distance metrics to construct distance vectors between pairs of persistence diagrams.  

\subsection{Vectorization and machine learning}
\begin{figure}[h]
    \centering
    \captionsetup{width=1\textwidth}
    \includegraphics[width=1\textwidth]{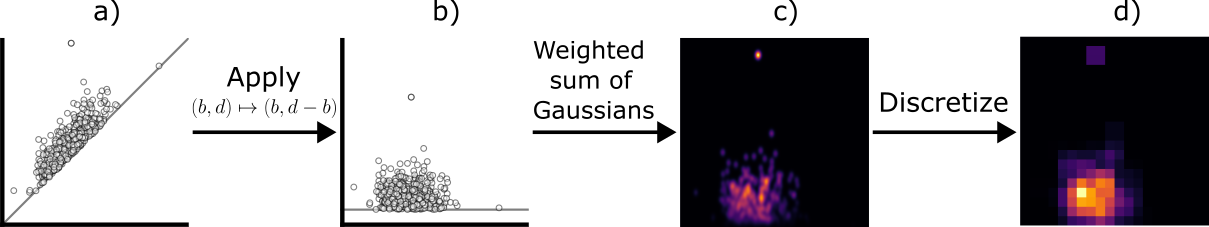}
    \caption{{\bf Vectorization of persistence diagrams via persistence image.}  a) An example persistence diagram. b) The result of mapping each point (birth, death) in a persistence diagram to (birth, death-birth). c) A weighted sum of Gaussians centered at each point of (b). d) A discretized array of image (c). The resulting persistence image is often flattened into a vector. Figure adapted from \cite{persistenceimage}.}
    \label{fig:PI_pipeline}
\end{figure}

Given a persistence diagram $\pd_k(X^\bullet)$, various techniques can be used to convert it into a vector that is compatible with standard statistics and machine learning \cite{TDA_vectors_survey, landscapes}. Here, we use persistence images \cite{persistenceimage}, which summarize the distribution of points on the persistence diagram using a weighted sum of Gaussian distributions centered at each point of the persistence diagram (see Fig. \ref{fig:PI_pipeline}). 

A persistence diagram is first transformed by mapping each point $(\text{birth}, \text{death})$ to \\$(\text{birth}, \text{death - birth})$ (Fig. \ref{fig:PI_pipeline}a,b). We then place a Gaussian distribution centered at each transformed point and assign a non-negative weighting function (Fig. \ref{fig:PI_pipeline}b,c). The function places zero weight for points along the horizontal axis of Fig. \ref{fig:PI_pipeline}b. The weighted sum of Gaussians is then discretized to produce an array called a persistence image (Fig. \ref{fig:PI_pipeline}d). The persistence image is often flattened into a vector. The resulting vector is influenced by several parameters, including the width, $\sigma$ of the Gaussian, and discretization size. In this study, we use $\sigma = 1$ and discretize images to size $20\times 20$, resulting in flattened vectors of dimension $400$.

\section{Introducing filtrations for multispecies data}
\label{sec:multisystem PH}

While standard PH detects structure in a point cloud, it fails to encode how multiple point clouds are related. We present two extensions of the standard PH pipeline to capture multi-system interactions: Dowker PH \cite{Chowdhury2018} and multispecies witness PH, a new construction motivated by witness complexes \cite{deSilva2004}. 
 
\subsection{Dowker persistent homology}
Let $U$ and $V$ denote two distinct point clouds. In our study, $U$ and $V$ represent different biological cell types, such as tumor cells and macrophages. The structure of $U$ from the viewpoint of $V$ can be studied using a Dowker filtration:

\begin{definition}[Dowker filtration \cite{dowker1952, Chowdhury2018}]
\label{def:Dowker complex}
Let $U$ and $V$ be point clouds, and let $d_{U,V}$ be the distance function between elements of $U$ and $V$. A \emph{Dowker complex at parameter $\varepsilon$}, denoted $\dowker_{U,V}^{\varepsilon}$, is a simplicial complex that has $U$ as the potential vertex set and includes the $n$-simplex $\sigma = (u_0,\ldots,u_n)$ if there exists a $v \in V$ such that $d_{U,V}(u_i,v) \leq \varepsilon$ for all $u_i \in \sigma$. A \emph{Dowker filtration} $\dowker_{U,V}^{\bullet}$ is a nested sequence of Dowker complexes $\dowker_{U,V}^{\varepsilon}$ for varying $\varepsilon$.
\end{definition}

\begin{figure}[h]
    \centering
    \captionsetup{width=1\textwidth}
    \includegraphics[width=0.9\textwidth]{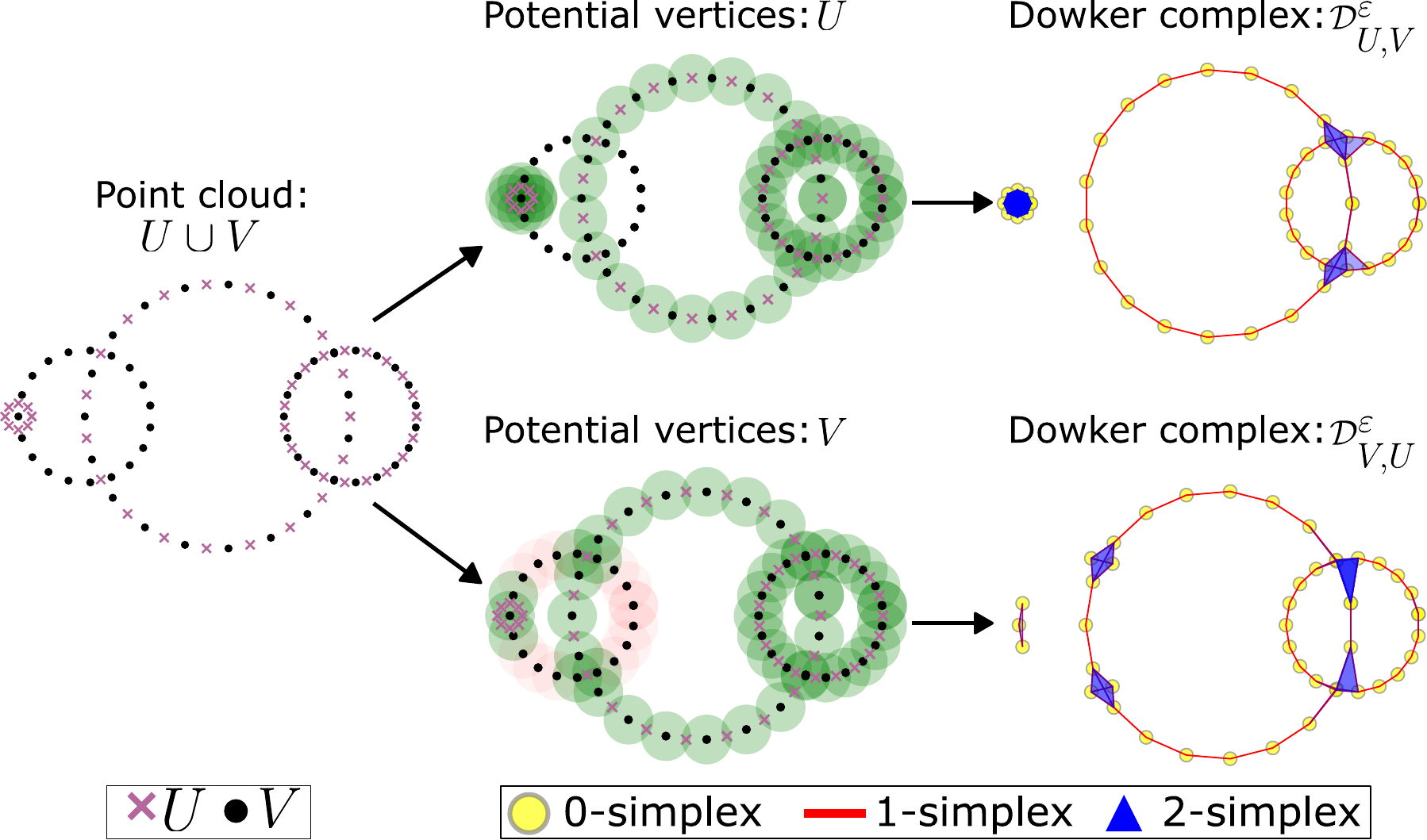}
    \caption{\textbf{Example Dowker complexes.} We present two Dowker complexes built on point clouds $U$ and $V$ for some proximity parameter $\varepsilon$. (Top) Dowker complex with $U$ as the potential vertex set. (Bottom) Dowker complex with $V$ as the potential vertex set. Given a potential vertex set, the $\varepsilon$-neighborhoods of the vertices are shown in green if the neighborhood contains an element of the other point cloud. Otherwise, the neighborhood is shown in red. A vertex with a green neighborhood becomes a 0-simplex in the Dowker complex. We add a 1-simplex between two vertices if their $\varepsilon$-neighborhood intersection contains a vertex from the other point cloud. We add a 2-simplex among three vertices if their $\varepsilon$-neighborhood intersection contains a vertex from the other point cloud.}
    \label{fig:dowker complex}
\end{figure}

The Dowker complex $\dowker_{U,V}^{\varepsilon}$ at parameter $\varepsilon$ captures relations between $U$ and $V$, where the relations are restricted to points $(u,v)$ whose distance is at most $\varepsilon$. Dowker complexes can capture shared topological features between two point clouds\footnote{There are instances in which the Dowker complex captures a feature present in $U$ that isn't present in $V$, for example, if $V$ is a dense sample of a region containing $U$. See Section \ref{sec:DPD_examples} for details.}, as illustrated in Fig. \ref{fig:dowker complex}. The Dowker complexes $\dowker_{U,V}^{\varepsilon}$ (Fig. \ref{fig:dowker complex}, top) and $\dowker_{V,U}^{\varepsilon}$ (Fig. \ref{fig:dowker complex}, bottom) each have $U$ and $V$ as the potential vertex set. Note that the two Dowker complexes resemble one another even though their vertex sets are distinct. For example, both Dowker complexes have two connected components and three 1-dimensional cycles, i.e., loops. Dowker's Theorem states that the two Dowker complexes have the same homology groups, i.e., connected components and loops \footnote{Note that while the homology groups of the above constructions are isomorphic, their connectivity, as measured for example by $Q$-analysis \cite{atkin1972cohomology}, may differ.} \cite{dowker1952}. In fact, the geometric realizations of the two Dowker complexes are homotopy equivalent \cite{BjornerTopMethods}.

To study the features of Dowker complexes across a range of parameters $\varepsilon$, we compute the PH of the Dowker filtration $\dowker_{U,V}^{\bullet}$. We call the resulting persistence diagram $\pd_k(\dowker_{U,V}^{\bullet})$ the Dowker persistence diagram. The functorial Dowker's Theorem states that the persistence diagrams of the two filtered Dowker complexes are the same. 

\begin{theorem}[Functorial Dowker's Theorem \cite{Chowdhury2018}]
\label{dowkertheorem}
$\pd_k(\dowker_{U,V}^{\bullet}) = \pd_k(\dowker_{V,U}^{\bullet})$ for all $k$.

\end{theorem}

The Dowker persistence diagram is a collection of birth and death parameters of $k$-dimensional topological features, i.e., connected components and loops for $k = 0$ and $k = 1$ respectively, in the Dowker filtration\footnote{We consider only $k = 0, 1$ in our analysis.}. The Dowker persistence diagram can be vectorized via persistence images as described in Section \ref{sec:single system TDA} and then be used in various statistical and machine learning methods.  

\subsection{Multispecies witness persistent homology}\label{Subsec:WPH}
Our second approach is motivated by the construction of \emph{(lazy) witness filtrations}.
The (lazy) witness filtration was first introduced by de Silva and Carlsson \cite{deSilva2004} and has been used to study noisy artificial datasets \cite{Kovacev2012}, primary visual cortex cell populations ~\cite{Singh2008}, and cancer gene expression data ~\cite{Lockwood2015}. Roughly, the lazy witness filtration is constructed via the following steps:
\begin{enumerate}
\item Select a subset of \textit{landmark points} $L$ from the point cloud $P$.
\item Construct a \textit{lazy witness filtration} where the landmarks $L$ are the vertex set and the full point cloud $P$ serve as \textit{witnesses} for higher order simplices. Broadly speaking, points in $P$ are witnesses to the simplices on $L$ to which they are closest.
De Silva and Carlsson \cite{deSilva2004} demonstrate that the resulting simplicial complex can be interpreted as an instrinsic Delaunay triangulation \cite{Delauney1934} of the point cloud. A filtration of the resulting simplicial complex is typically created by measuring the spatial scale of the simplices, similar to the Dowker filtration as described above\footnote{The Dowker filtration can be viewed as a special case of the lazy witness filtration. In the general formulation of the lazy witness filtration \cite{deSilva2004} the distance to the $\nu$-th closest witness is added to the proximity filtration scale $\epsilon$. Given point clouds $P$ and $Q$, a modified Dowker filtration in which all vertices have birth time $0$ is a witness filtration with $P$ as landmarks, $Q$ as witnesses, and $\nu = 0$. 

}. 
\end{enumerate}
For a multispecies point cloud $P = \cup_{i = 0}^m P_i$, $P_i \cap P_j = \emptyset$ for $i \neq j$, we use a similar construction to capture the spatial patterns of different $P_i$. However, rather than choosing a subset of landmarks $L$ from $P$, we use one of the point species as landmarks, i.e., $L = P_0$. Motivated by the close relationship of the witness complex and the Delaunay triangulation~\cite{deSilva2004}, we create the Delaunay triangulation \cite{Delauney1934} $D_0$ on the landmark set, i.e., for 2D point cloud data we create the triangulation of the 2D convex hull of $P_0$. We include all simplices from the Delaunay triangulation and their faces in our simplicial complex, i.e., for 2D data we include all triangles, their edges, and their vertices as the $2$-, $1$- and $0$-simplices of the simplicial complex. The remaining point species $P_i$ for $i = 1,...,m$ in $P$ are then used as witnesses for the simplices in the Delaunay triangulation:

\begin{definition}[$P_i$-witness point]
Let $p \in P_i$, $l \in L$, and $d$ a distance function on $P$. We say that $p$ is a \emph{$P_i$-witness} for 
the $n$-simplex $\sigma = (l_0,\dots,l_n)$ if $d(p,l_i) \leq d(p,\hat{l})$ for all $\hat{l} \in L\setminus\{l_0,\dots,l_n\}$ and $i = 0,...,n$.
\end{definition}

We now create species-dependent filtrations $W^\bullet_{0,i}$ on the landmark set $P_0$ using witness points from $P_i$:

\begin{definition}[Multispecies witness filtration $W^\bullet_{0, i}$]
Let $P = \cup_{i = 0}^m P_i$ denote a collection of different point clouds, and let $D_0$ be the Delaunay triangulation of $P_0$. The \emph{multispecies witness filtration} is a sequence of nested simplicial complexes $W^\bullet_{0, i}$ on $P_0$ with respect to witness points in $P_i$ where $W^\mu_{0,i}$ has $P_0$ as its potential vertex set and includes the $n$-simplex $\sigma = (p_0,\dots,p_n) \in D_0$ and all its faces, if $\tilde{\mu}^\sigma \leq \mu$ with $\tilde{\mu}^\sigma = \frac{\mu_{\text{max}}-\mu^\sigma}{\mu_{\text{max}}}$, where $\mu^\sigma$ is the number of $P_i$-witnesses of $\sigma$ and $\mu_{\text{max}}$ is the maximal number of $P_i$-witnesses for a simplex in $D_0$.
\end{definition}

\begin{figure}[h!]
\centering
\captionsetup{width=1\textwidth}
\includegraphics[width=1\textwidth]{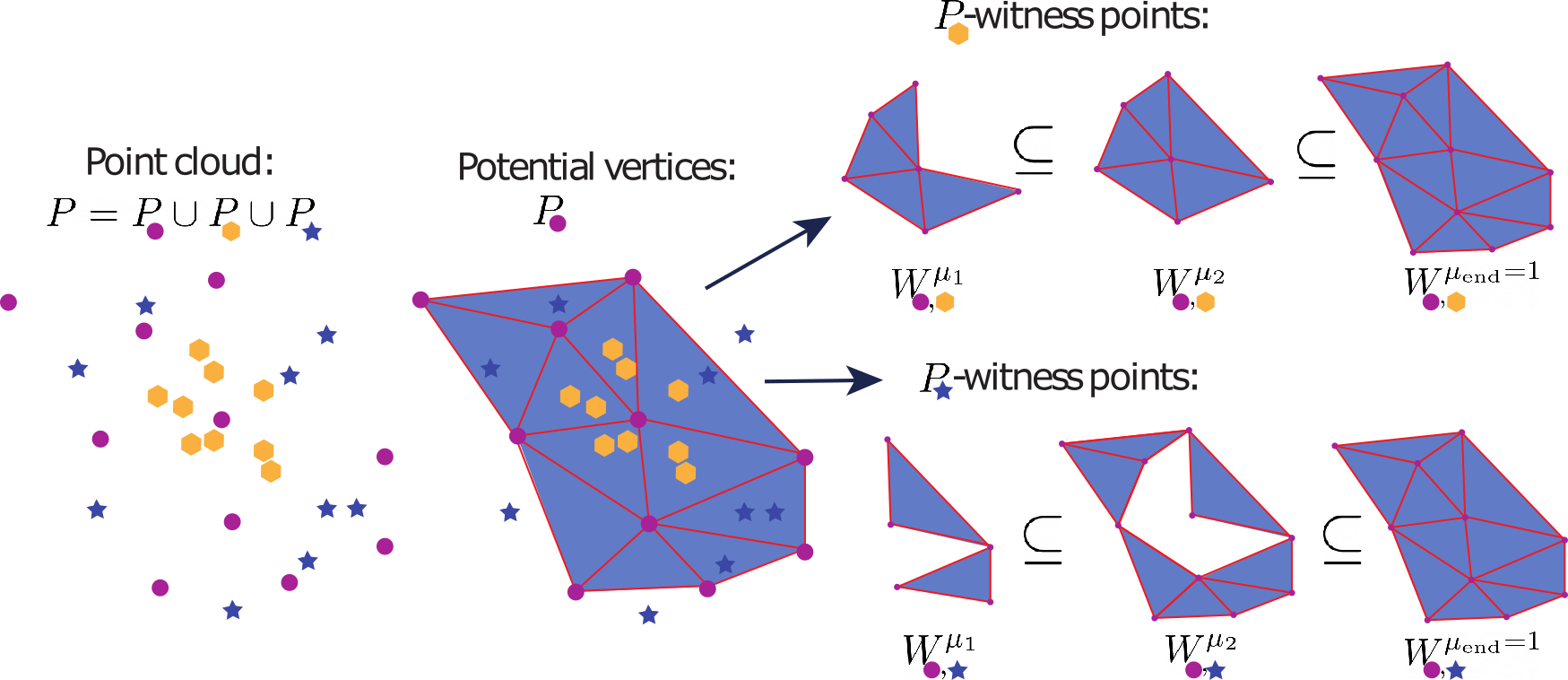}
    \caption{{\bf Example multispecies witness filtration.} Given a point cloud $P = P_0 \cup P_1 \cup P_2$, we illustrate two multispecies witness filtrations on the Delaunay triangulation on $P_0$ using witness points from $P_1$ (depicted as yellow hexagons) and witness points from $P_2$ (depicted as blue stars). Different witness points give rise to different filtrations of the same simplicial complex. 
    \label{fig:witnessPH}}
\end{figure}

We illustrate the multispecies witness filtration in an example point cloud in Fig. \ref{fig:witnessPH}. 

To compare the effect of the different types of witnesses on the filtration, we first compute the dimension-0 and dimension-1 persistence diagrams of the multispecies witness filtrations, denoted $\pd_0(W^\bullet_{0,i})$ and $\pd_1(W^\bullet_{0,i})$, for $i=1, \dots, m$, and we compute pairwise distance vectors among the different persistence diagrams. We focus on distances between persistence diagrams. The entries of our distance vectors are given by the pairwise Bottleneck distances $d_B$ among $\pd_0(W^\bullet_{0,i})$, the pairwise Bottleneck distances $d_B$ among $\pd_1(W^\bullet_{0,i})$, the pairwise 1-Wasserstein distances $d_W$ among $\pd_0(W^\bullet_{0,i})$, and the pairwise 1-Wasserstein distances $d_W$ among $\pd_1(W^\bullet_{0,i})$ for $i=1, \dots, m$. Given a point cloud $P = \bigcup_{i=0}^m P_i$ with $m+1$ species, this results in distance vectors with $2 \times 2 \times {m \choose 2} = 2m(m-1)$ entries. Remark, this choice of distance vector sidesteps the additional steps (and parameter choices) of constructing persistence image-based distances because different witness points lead to differences manifesting in the filtrations of the Delaunay triangulation of $P_0$ (see Fig. \ref{fig:witnessPH}).

\section{Results}
\label{sec:results}
We demonstrate the utility of relational PH in predicting the macrophage phenotype (Problem 1) and in classifying the qualitative behavior of different parameter regimes of the ABM (Problem 2). For the first task, we find that using Dowker PH features improves the performance of a classifier in comparison to using both non-relational topological and non-topological features. In particular, we find that Dowker PH between tumor cells and blood vessels is the best predictor for the dominant macrophage phenotype. For the second task, we perform classification using the multispecies witness filtration features and recover the previous subjective classification of Fig.~\ref{fig:joshlabelsII}. 

\subsection{Dowker persistent homology predicts dominant macrophage phenotype}
\label{sec:results_Dowker}
\subsubsection{Prediction pipeline}

We classify a synthetic tumor microenvironment as either anti-tumor ($M_1$) macrophage dominant or pro-tumor ($M_2$) macrophage dominant based on the spatial distributions of blood vessels, tumor cells, and macrophages. Since the $M_1$ and $M_2$ macrophages exhibit significantly different dynamics in the tumor microenvironment (see Section \ref{sec:ABM}), we hypothesize that the relations of spatial distributions among the three cell types are good predictors of the dominant macrophage phenotype. Our input data is a point cloud $P = P_V \cup P_T \cup P_M$ that represents the locations of the three cell types. Note that the input data is blind to the phenotype of individual macrophages. 

Given a point cloud $P$, if 50\% or more macrophages are $M_1$ macrophages, then we label the point cloud as $M_1$ dominant. Otherwise, we label the point cloud as $M_2$ dominant. A total of 731 images are labeled 0 ($M_1$ dominant), and 241 images are labeled 1 ($M_2$ dominant). 

\begin{figure}[H]
    \centering
    \captionsetup{width=1\textwidth}
    \includegraphics[width=\textwidth]{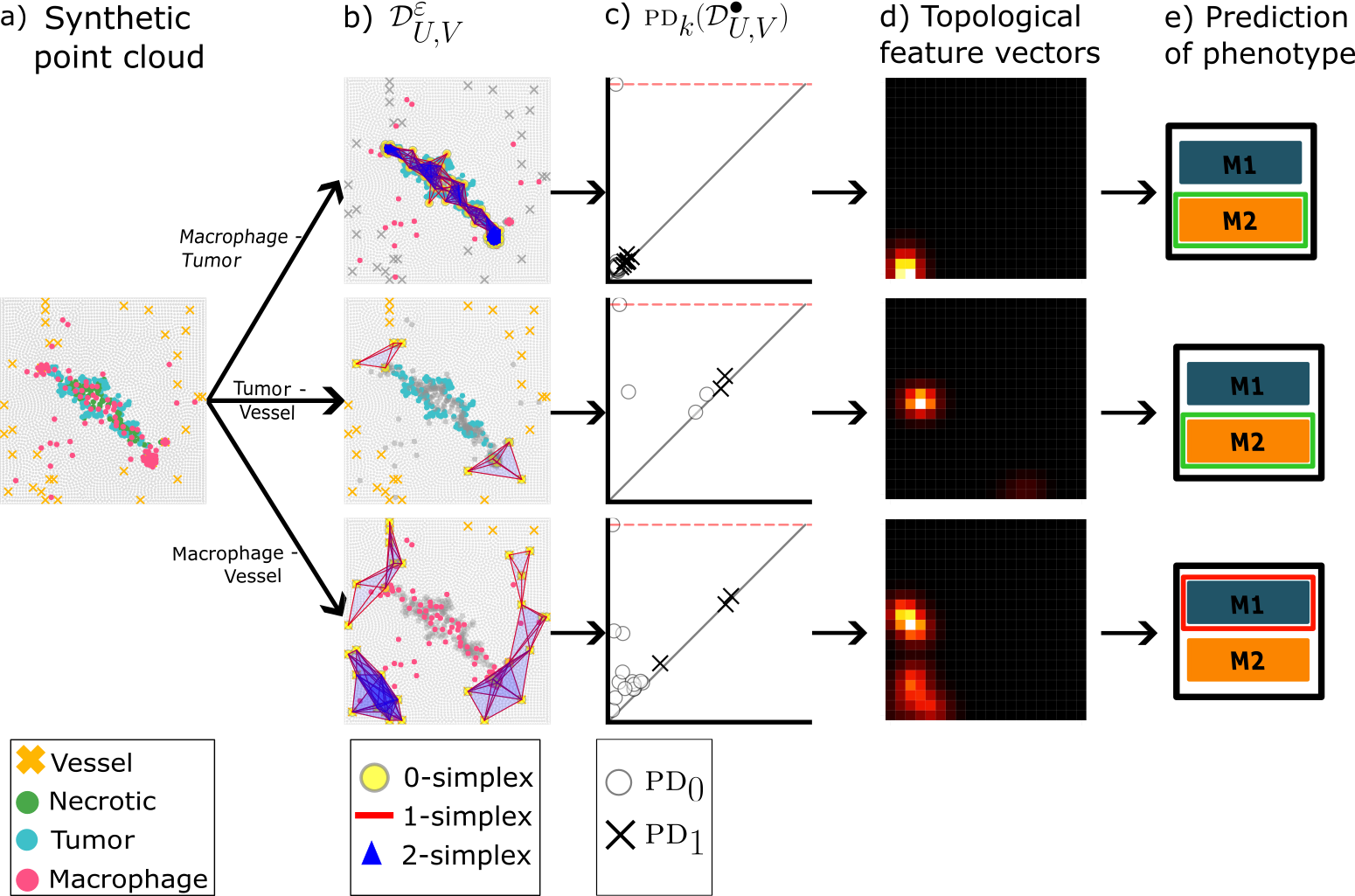}
    \caption{\textbf{Pipeline for macrophage phenotype prediction using Dowker PH.} a) A point cloud representing a synthetic tumor microenvironment generated by an ABM. b) Dowker complexes built on different pairs of cells at fixed proximity parameters. c) Dowker persistence diagrams $\pd_k(\dowker_{U,V}^{\bullet})$ for $k = 0, 1$. d) Vectorization of (dimension-0) Dowker persistence diagrams via persistence images. e) An SVM classifier takes a flattened persistence image as input and predicts the dominant macrophage phenotype of the synthetic tumor microenvironment.}
    \label{fig:dowker pipeline}
\end{figure}

For each $P$, we use Dowker PH to capture relations between pairs of constituents of the tumor microenvironment\footnote{We only used spatial information of macrophages to create topological descriptors. Macrophage phenotype information is used only to label images.} (see Fig.~\ref{fig:dowker pipeline}a,b). We consider the following three pairs of cell types: macrophages and tumor cells, tumor cells and blood vessels, and macrophages and blood vessels (see Fig. \ref{fig:dowker pipeline}b). For each pair, we compute the dimension-0 and dimension-1 Dowker persistence diagrams\footnote{Recall that the Dowker persistence diagram is agnostic to the choice of the vertex set (Theorem \ref{dowkertheorem}). In each pair, we chose the cell type with a smaller number of points as the vertex set for faster computation.} (see Fig. \ref{fig:dowker pipeline}c). Each point cloud thus results in six Dowker persistence diagrams: $\pd_0(\dowker_{M,V}^{\bullet})$, $\pd_1(\dowker_{M,V}^{\bullet})$, $\pd_0(\dowker_{T,V}^{\bullet})$, $\pd_1(\dowker_{T,V}^{\bullet})$, $\pd_0(\dowker_{M,T}^{\bullet})$, $\pd_1(\dowker_{M,T}^{\bullet})$.

Each Dowker persistence diagram is vectorized via persistence images to an array of size 20 × 20 \footnote{In this study, the classification accuracy is fairly robust to the size of the persistence image. Such robustness is known in the literature \cite{persistenceimage}.} (see Fig.~\ref{fig:dowker pipeline}d). We flatten the persistence images into vectors of size $400$ and train a Support Vector Machine (SVM) for the image classification task. (see Fig.~\ref{fig:dowker pipeline}e).

We also train SVMs on non-relational topological features obtained from four Vietoris-Rips persistence diagrams: $\pd_0(\VR^{\bullet}_{T})$, $\pd_1(\VR^{\bullet}_{T})$, $\pd_0(\VR^{\bullet}_{M})$, $\pd_1(\VR^{\bullet}_{M})$.
We further train an SVM on non-topological features such as the count of each cell type and the average distance of each cell type to the nearest blood vessels (see data description in Section \ref{sec:ABM}). 

For each SVM classifier, we optimize the hyperparameters via stratified 5-fold cross-validation, employing the synthetic minority oversampling technique (SMOTE) \cite{smote} in each fold to address the class imbalance. We train an SVM on 10 different random splits of train and test data and report the 10 classification accuracies on the test data. 

\subsubsection{Dowker persistence diagrams capture shared topological features}
\label{sec:DPD_examples}

\begin{figure}[h]
    \centering
    \captionsetup{width=1\textwidth}
    \includegraphics[width=\textwidth]{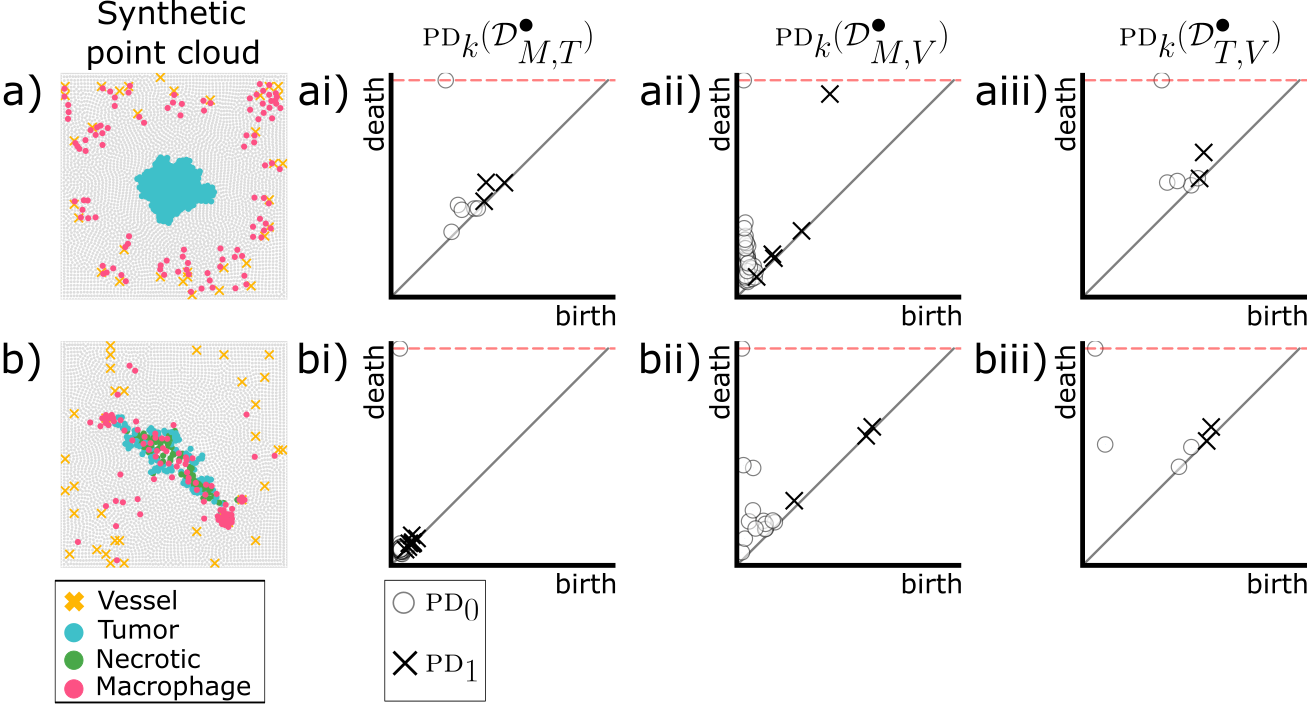}
    \caption{\textbf{Dowker persistence diagrams capture spatial relations between cell types.} a) A synthetic tumor microenvironment in which macrophages and blood vessels surround a compact tumor. The six Dowker persistence diagrams are shown in ai, aii, aiii. b) A synthetic tumor microenvironment where the cancer cells and macrophages occupy different spaces from the blood vessels. The cancer cells and macrophages are in close proximity to blood vessels in two regions -- the top left and bottom right corners of the tumor mass. The six Dowker persistence diagrams are shown in bi, bii, biii. ai) The large birth parameters of points in $\pd_0(\dowker_{M,T}^\bullet)$ indicate that macrophages and tumor cells are far from one another. aii) The small birth parameters of points in $\pd_0(\dowker_{M,V}^{\bullet})$ indicate that macrophages and blood vessels are colocalized. The single cross far from the diagonal in $\pd_1(\dowker_{M,V}^{\bullet})$ indicates that macrophages and blood vessels share a common loop. aiii) Both Dowker persistence diagrams are similar to the diagrams in panel (ai) because the relationship between tumor cells and blood vessels is similar to the relationship between macrophages and tumor cells. bi) The small birth parameters of $\pd_0(\dowker_{M,T}^{\bullet})$ indicate that macrophages and tumor cells occupy similar regions. bii) The spread of birth parameters for points in $\pd_0(\dowker_{M,V}^{\bullet})$ indicates the variance in the extent to which macrophages and vessels occupy similar spaces. biii) The two points in $\pd_0(\dowker_{T,V}^{\bullet})$ far from the diagonal indicate that there are two regions (the top left and bottom right corners of the tumor mass) where the tumor cells and the blood vessels are close to each other.}
    \label{fig:dowker examples}
\end{figure}

Before we discuss classification accuracy, we present example point clouds and interpretation of Dowker persistence diagrams (see Fig. \ref{fig:dowker examples}).

Recall that $\pd_0(\dowker_{U,V}^{\bullet})$ summarizes the birth and death of connected components of Dowker complexes as one varies the distances between $P_U$ and $P_V$. One can thus consider a dimension-0 Dowker persistence diagram as summarizing shared connected components between two point clouds. There are multiple ways in which a shared connected component arises - $P_U$ and $P_V$ might occupy a similar region, or $P_U$ and $P_V$ may occupy different regions but have close contact. In such cases, the shared features will be represented by points in $\pd_0(\dowker_{U,V}^{\bullet})$ with small birth parameters. 

For example, consider the relationship between macrophages and tumor cells in Fig.~\ref{fig:dowker examples}a and Fig.~\ref{fig:dowker examples}b. In Fig.~\ref{fig:dowker examples}a, the macrophages are distant from the tumor cells, so the points in $\pd_0(\dowker_{M,T}^{\bullet})$ have large birth times (see Fig. \ref{fig:dowker examples}ai). On the other hand, in Fig.~\ref{fig:dowker examples}b, the macrophages and tumor cells occupy similar spaces, so the points in $\pd_0(\dowker_{M,T}^{\bullet})$ have small birth times (see Fig. \ref{fig:dowker examples}bi). 

Consider the relationship between tumor cells and blood vessels in Fig. \ref{fig:dowker examples}b. The tumor cells and blood vessels mostly occupy different spaces. However, the tumor cells and blood vessels are in close proximity in two regions, one on the top left corner and another on the bottom right corner of the tumor mass. The fact that there are two ``contact points" between the tumor and blood vessels is reflected by two points in $\pd_0(\dowker_{T,V}^{\bullet})$ that are far from the diagonal (Fig.~\ref{fig:dowker examples}biii). In Fig. \ref{fig:dowker examples}a, the macrophages and blood vessels occupy very similar regions. Such colocalization between macrophages and blood vessels is reflected by the abundance of points in $\pd_0(\dowker_{M,V}^{\bullet})$ with small birth times (see Fig. \ref{fig:dowker examples}aii).

A dimension-1 Dowker persistence diagram summarizes the evolution of cycles of Dowker complexes as one varies the distances between $U$ and $V$. We interpret points in $\pd_1(\dowker_{U,V}^{\bullet})$ that are far from the diagonal line as representing shared loops between two point clouds
\footnote{We caution the reader that $\pd_1(\dowker_{U,V}^{\bullet})$ can contain points far from the diagonal line even if $P_U$ and $P_V$ do not necessarily have shared cycles. Such a situation arises, for example, when $P_U$ is sampled from a circle while $P_V$ is a dense, uniform sample of the background.}. For example, the macrophages and blood vessels in Fig.~\ref{fig:dowker examples}a share a loop structure, and such shared loop is reflected by a point in $\pd_1(\dowker_{M,V}^{\bullet})$ that is far from the diagonal (Fig. \ref{fig:dowker examples}aii).

\subsubsection{SVM on Dowker features predicts dominant macrophage phenotype}
\begin{figure}[H]
    \centering
    \captionsetup{width=1\textwidth}
        \includegraphics[width=\textwidth]{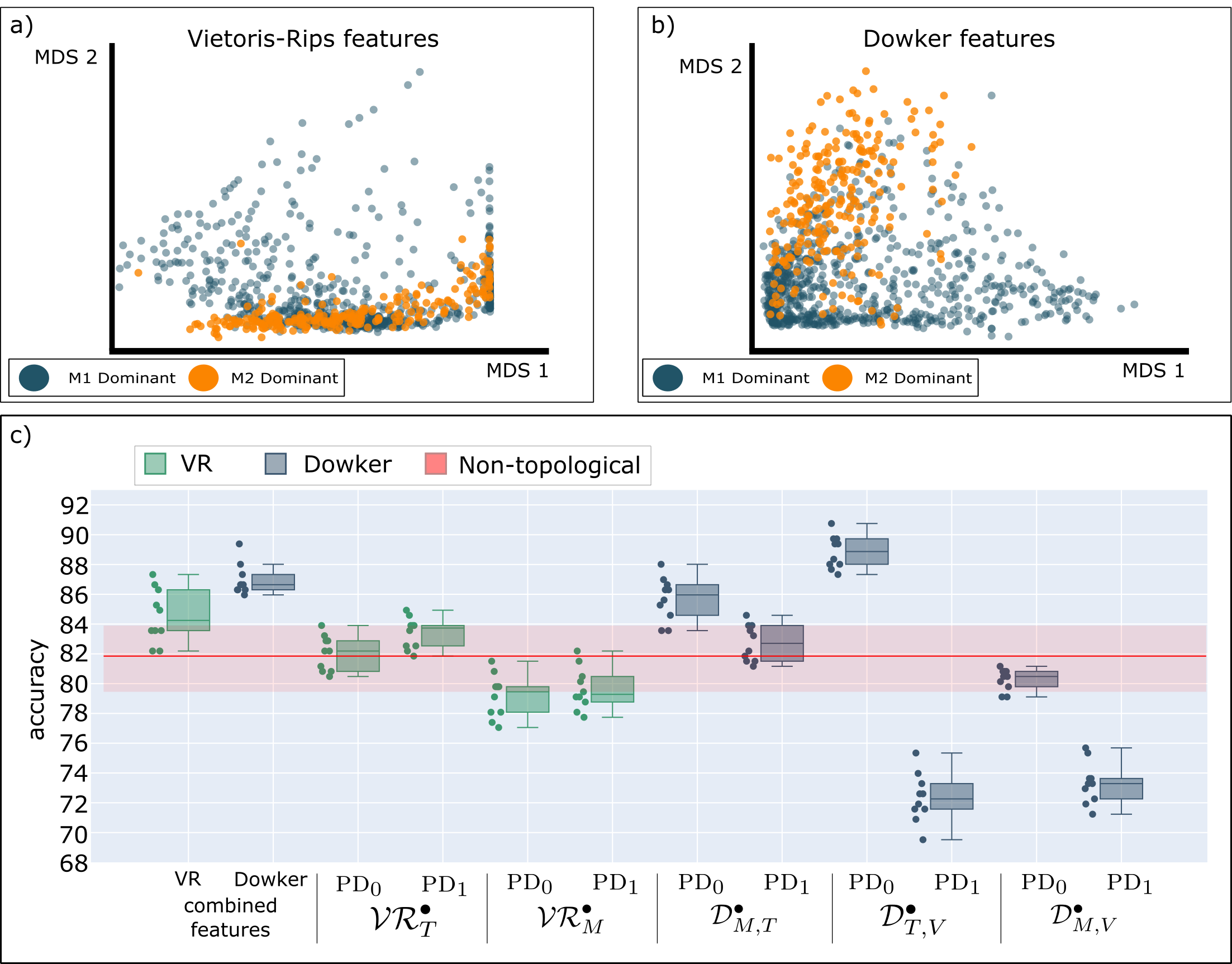}
    \caption{{\bf Dowker persistent homology features improve the 
   prediction of dominant macrophage subtype.} a) MDS projection of Vietoris-Rips features. b) MDS projection of Dowker features. The two classes have better separation when using Dowker features than the Vietoris-Rips features. c)  Classification accuracies of SVMs trained on Vietoris-Rips features (green), Dowker features (navy), and non-topological features (red). The box plot summarizes the accuracies from 10 different splits of train and test data. The red box shows the minimum (lower bounding line), median (middle line), and maximum (upper bounding line) accuracy values for SVM trained on non-topological feature vectors. The first two box plots show the accuracies of two SVMs, one trained on all Vietoris-Rips features and another trained on all Dowker features. SVM trained on Dowker features has higher accuracy than SVM trained on Vietoris-Rips features. The remaining box plots show the accuracies of SVMs trained on individual Vietoris-Rips or Dowker features. SVM trained on $\pd_0(\dowker^\bullet_{T, V})$ has the highest accuracy among all SVM trained on the Dowker features. }
    \label{fig:svm_accs}
\end{figure}

We first visually inspected whether Dowker persistence diagrams can distinguish $M_1$ and $M_2$ dominant tumor microenvironments. Recall that we computed six Dowker persistence diagrams, which resulted in six $400$-dimensional vectors. We concatenated the six vectors into a $2400$-dimensional vector, and we refer to the resulting vector as a Dowker feature vector. A two-dimensional visualization via Multidimensional Scaling (MDS) \cite{Kruskal1964MultidimensionalSB} shows decent separation of classes (see Fig. \ref{fig:svm_accs}b). For comparison, we computed four Vietoris-Rips persistence diagrams from tumor cells and macrophages, vectorized, and concatenated vectors. We refer to the concatenated vectors as Vietoris-Rips features. A comparison of MDS on the Vietoris-Rips features (Fig. \ref{fig:svm_accs}a) indicates that Dowker features may be better predictors of the dominant macrophage phenotype.

We train two SVM classifiers, one that takes the Dowker feature vectors as input and another that takes the Vietoris-Rips feature vectors as input. The SVM trained on Dowker features has higher accuracy (median accuracy $86.6\%$) than the SVM trained on Vietoris-Rips features (median accuracy $84.2\%$). Furthermore, the lower quartile of accuracy from Dowker features is roughly equal to the upper quartile of accuracy from Vietoris-Rips features ($\sim 86\%$) (see Fig. \ref{fig:svm_accs}c). Both models outperform an SVM trained on non-topological features such as the number of cells per cell type and average distances of cell types to the nearest blood vessels (see Fig. \ref{fig:svm_accs}c).

Next, we investigate which cell types were most informative in predicting the dominant macrophage of the synthetic tumor microenvironment. To this end, we trained ten additional SVM classifiers. We train four classifiers on the four Vietoris-Rips features and six classifiers on the six Dowker features. Among the classifiers trained on Vietoris-Rips features, the model trained on $\pd_1(\VR^{\bullet}_{T})$ has the highest median accuracy ($83.7 \%$). One possible explanation is that $M_2$ macrophages assist metastasis of tumor cells by guiding them away from the tumor mass towards the blood vessels. During this process, the tumor cells may create many small loops as they navigate away from the tumor mass, creating many non-trivial points in $\pd_1(\VR^{\bullet}_{T})$. The persistence diagram $\pd_1(\VR^{\bullet}_{T})$ may then reflect the extent to which $M_2$ macrophages assist the spread of cancer cells. 

Among the classifiers trained on Dowker features, the model trained on $\pd_0(\dowker_{T,V}^{\bullet})$ has the highest accuracy (median accuracy $88.9\%$), followed by the model trained on $\pd_0(\dowker_{T,M}^{\bullet})$ ($86.0 \%$). It is perhaps surprising that the best predictor of the dominant macrophage phenotype does not involve the spatial distribution of macrophages. One possible explanation for the improved performance of models using $\pd_0(\dowker_{T,V}^{\bullet})$ is that $\pd_0(\dowker_{T,V}^{\bullet})$ indicates colocalization between tumor cells and blood vessels, which can represent the extent to which $M_2$ macrophages have assisted the tumor cells to navigate towards blood vessels for metastasis. 

Note that dimension-1 Dowker features involving blood vessels are not particularly good predictors of the dominant macrophage phenotype (see Fig.~\ref{fig:svm_accs}c). The poor performance may be due to the lack of common loops between blood vessels and tumor cells and between the blood vessels and the macrophages.

\newpage
\subsection{Multispecies witness features identify qualitative model behaviors}
\label{sec:WitnessFeatures}

\begin{figure}[H]
    \centering
    \captionsetup{width=1\textwidth}
        \includegraphics[width=\textwidth]{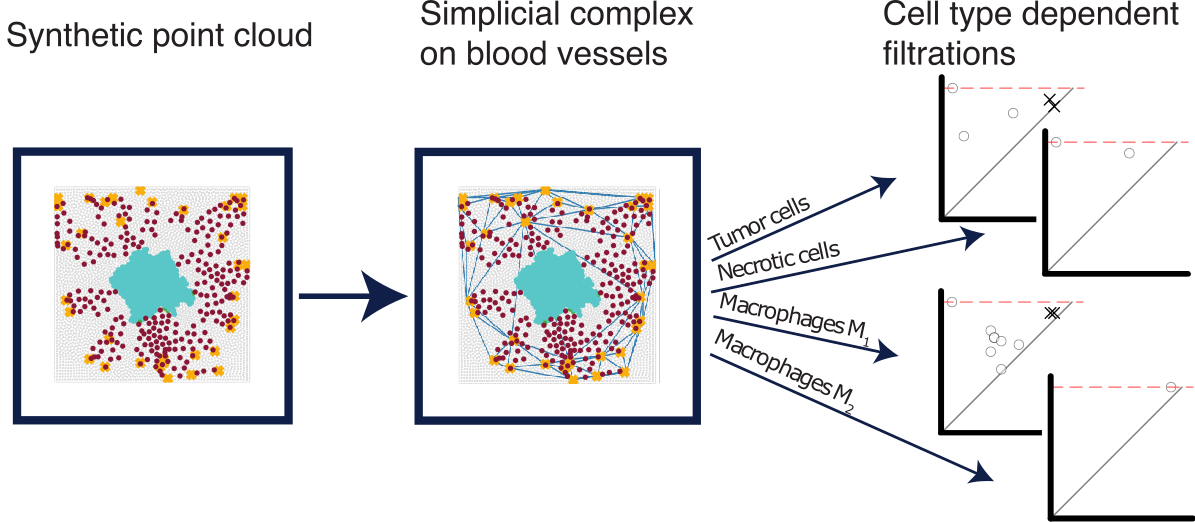}
    \caption{{\bf Multispecies witness PH on synthetic data from ABM.} The point cloud given by the synthetic data $P = P_{\text{V}} \cup P_{\text{T}} \cup P_{\text{N}} \cup P_{M_1} \cup P_{M_2}$ consists of blood vessels $P_V$, tumor cells $P_{\text{T}}$, necrotic cells $P_{\text{N}}$, anti-tumor macrophages $P_{M_1}$ and pro-tumor macrophages $P_{M_2}$. We construct a Delaunay triangulation on the blood vessels $P_{\text{V}}$ and build cell type dependent filtrations $W_{V,i}^\bullet$ of the Delaunay triangulation where $i \in \{\text{T},\text{N},M_1,M_2 \}$. We obtain one persistence diagram for each cell type specific filtration. 
    \label{fig:pipelinewitnessPH}}
\end{figure}

To study the different qualitative behaviors of the ABM, 
we focused on differences between the spatial distributions of the different cell types and applied the multispecies witness PH. We illustrate how we applied multispecies witness PH to the output of our ABM in Fig. \ref{fig:pipelinewitnessPH}. 

\begin{figure}[H]
    \centering
    \captionsetup{width=1\textwidth}
        \includegraphics[width=.6\textwidth]{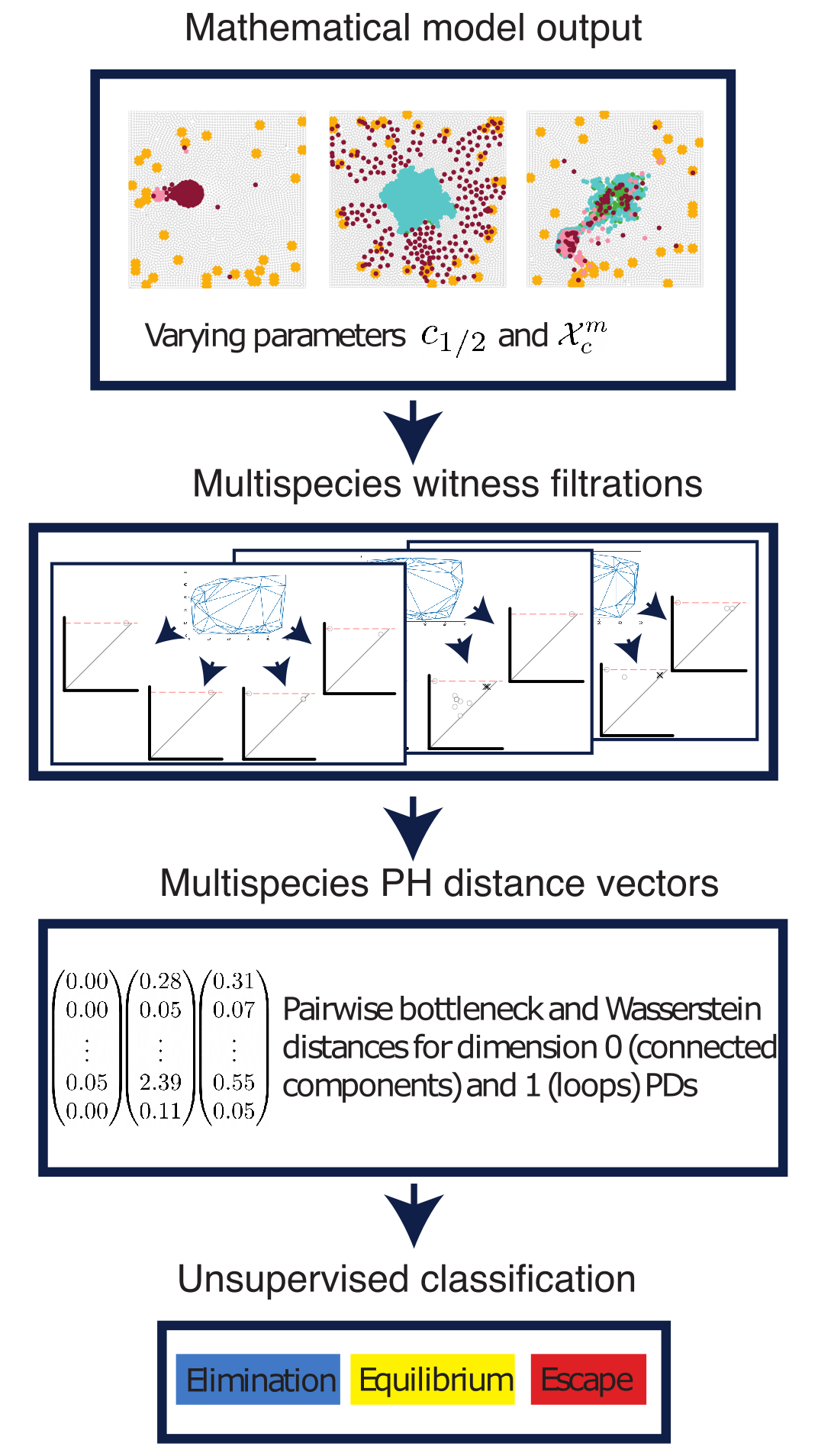}
    \caption{{\bf Multispecies witness PH pipeline.} We use the point cloud generated by an ABM as input into our multispecies witness filtrations. We compute persistence diagrams for the multispecies witness filtrations, thereby obtaining topological descriptors of the spatial heterogeneity in the input images. We use the persistence diagrams to compute multispecies PH distance vectors. The entires of these vectors correspond to the pairwise Bottleneck and 1-Wasserstein distances between the dimension-0 and dimension-1 persistence diagrams of the cell type specific filtrations. We use the multispecies PH distance vectors as input into unsupervised classification to identify different qualitative behaviors of the ABM.}
    \label{fig:pipeline}
\end{figure}

Our point cloud data $P = P_{\text{V}} \cup P_{\text{T}} \cup P_{\text{N}} \cup P_{M_1} \cup P_{M_2}$ consists of blood vessels $P_V$, tumor cells $P_{\text{T}}$, necrotic cells $P_{\text{N}}$, anti-tumor macrophages $P_{M_1}$ and pro-tumor macrophages $P_{M_2}$. We chose to fix $P_0 = P_V$ and considered two different versions for the witness filtrations: first, we did not distinguish macrophage phenotype, i.e., all macrophages are assumed to be identical and $P_{\text{M}} = P_{M_1} \cup P_{M_2}$. We obtained three different witness filtrations using tumor cells, necrotic cells, and macrophages as witness points. In the second case, we distinguished $M_1$ and $M_2$ macrophage subtypes and constructed four witness filtrations using tumor cells, necrotic cells, $M_1$ macrophages, and $M_2$ and macrophages as witness points. From the persistence diagrams, we computed multispecies PH distance vectors (see Subsection \ref{Subsec:WPH}) to compare the effect of the different types of witnesses on the filtration. 
The entries of our distance vectors are listed in Table \ref{tab:versions}. The pairwise distances each contributed 3 entries when all macrophages are considered to be the same cell type and 6 entries when distinguishing between $M_1$ and $M_2$ macrophages for each topological dimension considered. In this way, we converted each point cloud $P$ into a 12- (version 1) and a 24-dimensional (version 2) distance vector, respectively (for a summary, see Table \ref{tab:versions}). We used these distance vectors as input into $k$-means clustering. We summarize the full multispecies witness PH pipeline in Fig. \ref{fig:pipeline}. We compared our results to clustering performed on simple (non-topological) descriptor vectors (see data description in Section \ref{sec:ABM} for description of simple vectors and see Fig. \ref{fig:trivial clusters} in the Appendix for results).

\subsubsection{Multispecies witness persistence classification disregarding macrophage subtype}

\begin{figure}[H]
    \centering
    \captionsetup{width=1\textwidth}
        \includegraphics[width=\textwidth]{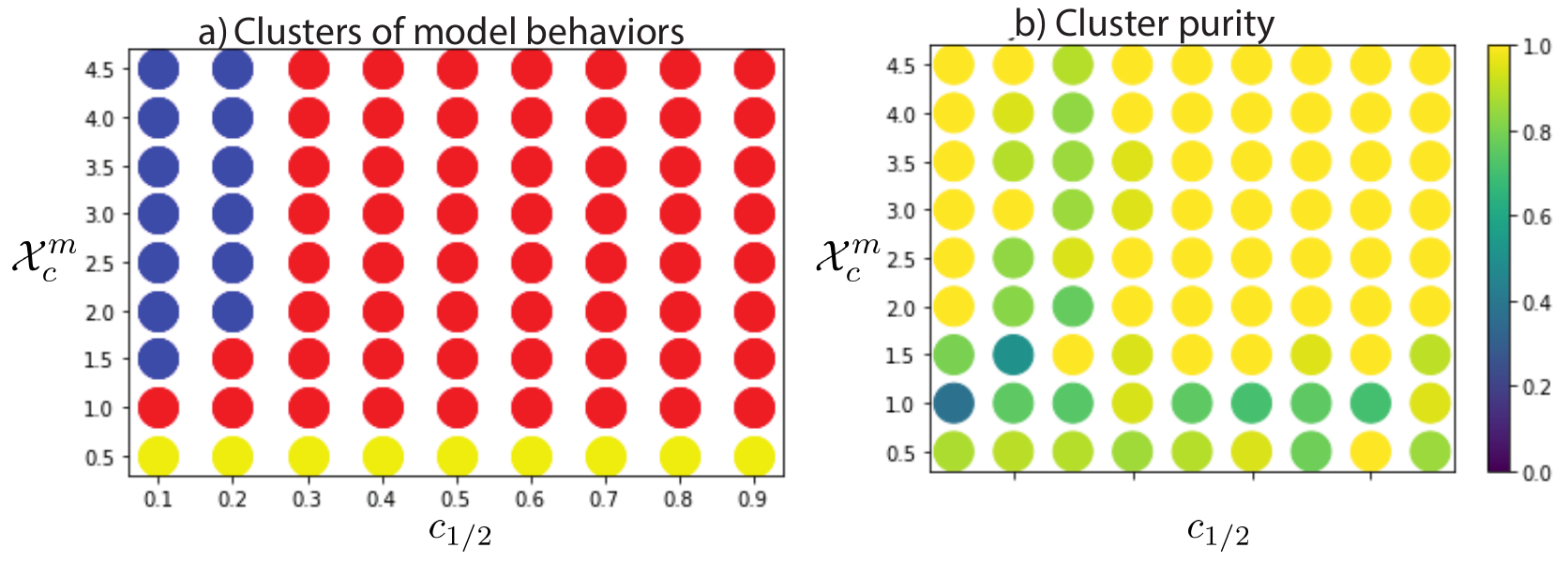}
    \caption{{\bf Classification of multispecies PH distance without distinguishing between macrophage subtypes.} a) Classification results. b) Cluster purity scores. For each parameter combination $\mathcal{X}^m_c$ and $c_{1/2}$ of the ABM, we include 20 independent simulations in our analysis. The colors red, blue, and yellow represent the cluster to which the majority of simulations are attributed by the $k$-means algorithm for $k=3$. The purity score is computed by taking the ratio between the number of simulations attributed to the majority clusters by 20. 
    \label{fig:MPHresults}}
\end{figure}

We recovered the three qualitatively different behaviors of the ABM using the unsupervised multispecies witness PH pipeline without including knowledge about macrophage subtypes. We applyed $k$-means classification for $k = 3$. Fig. \ref{fig:MPHresults} shows which of the three clusters is dominant amongst the 20 simulations for each parameter combination of $\mathcal{X}^m_c$ and $c_{1/2}$ that we consider. The results are consistent with the subjective classification of the qualitative behaviors of the model shown in Fig. \ref{fig:joshlabelsII}, i.e., we recovered parameter regimes dominated by tumor elimination, tumor macrophage equilibrium, and escape of the tumor, with the exception of simulations in regimes at the boundaries between the three behaviors. We investigated the consistency of the cluster assignment, which we refer to as \emph{cluster purity} by dividing the number of simulations attributed to the majority cluster by the total number of simulations for the parameter combination.
We found that cluster assignment is less consistent in simulations of the ABM that lie in boundary regions between different qualitative behaviors than in parameter regimes far away from boundaries (see Fig. \ref{fig:MPHresults}). Our results clearly surpass clustering obtained using simple descriptor vectors of the data (see Fig. \ref{fig:trivial clusters} in the Appendix), including information such as the number of cells per cell type and average distances of cell types to the nearest blood vessels with respect to cluster consistency with the subjective clusters shown in Fig. \ref{fig:joshlabelsII}.

\subsubsection{Multispecies witness persistence classification including macrophage subtypes}

\begin{figure}[H]
    \centering
    \captionsetup{width=1\textwidth}
        \includegraphics[width=\textwidth]{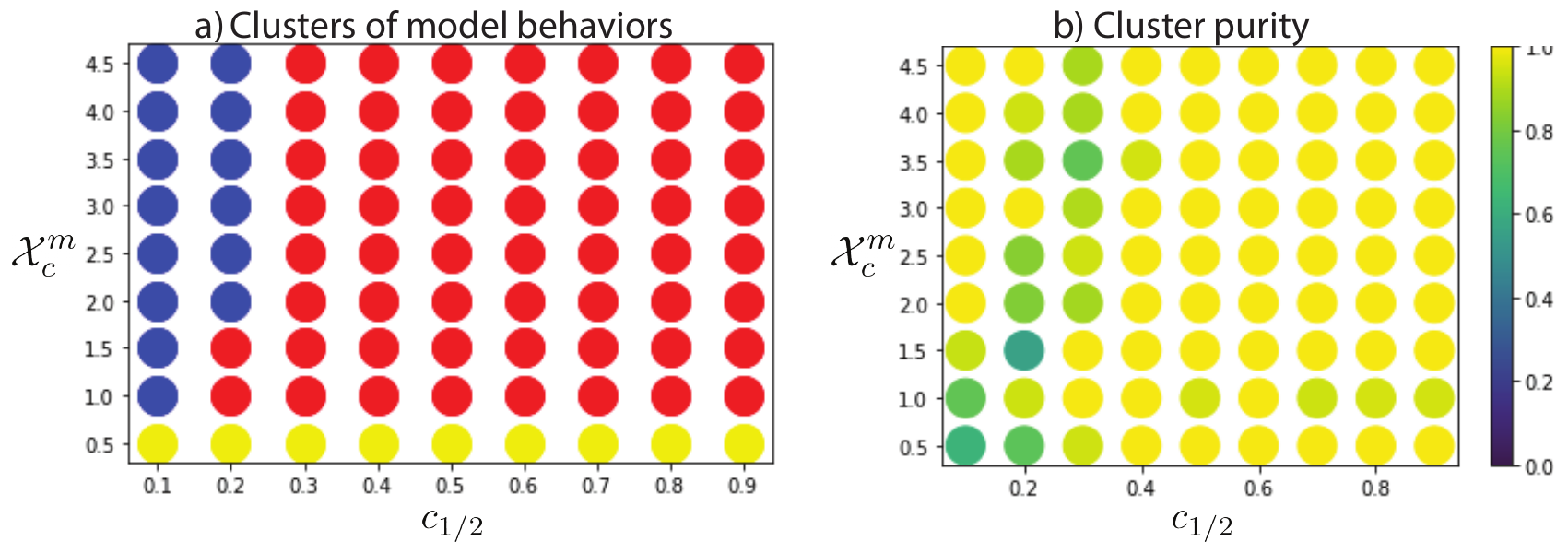}
    \caption{{\bf Classification of multispecies PH distance vectors distinguishing between macrophage subtypes $M_1$ and $M_2$.} a) Classification results. b) Cluster purity scores. For each parameter combination $\mathcal{X}^m_c$ and $c_{1/2}$ of the ABM, we include 20 independent simulations in our analysis. The colors red, blue, and yellow represent the cluster to which the majority of simulations are attributed by the $k$-means algorithm for $k=3$. The purity score is computed by taking the ratio between the number of simulations attributed to the majority clusters by 20.
    \label{fig:MPHresultsM1M2}}
\end{figure}

We also recovered the three qualitatively different behaviors of the ABM when information about macrophage subtypes $M_1$ and $M_2$ is included in the construction of our multispecies PH distance vectors. We show our results in Fig. \ref{fig:MPHresultsM1M2}. Comparison of the results in Fig. \ref{fig:MPHresults} and Fig. \ref{fig:MPHresultsM1M2} shows that the inclusion of the additional information about macrophage subtype alters the prediction of the qualitative behaviors for only one parameter combination, $\mathcal{X}^m_c = 1$ and $c_{1/2} = 0.1$, which is located at the phase transition between elimination and escape. We also computed the purity of clusters for each parameter combination by dividing the number of simulations attributed to the majority cluster by the total number of simulations for the parameter combination. We find that clusters assigned to parameter combinations located at the phase transitions between different parameter regimes are less consistent than those far away from boundaries. Again, our results surpass clustering obtained using simple descriptor vectors of the data, including information such as the number of cells per cell type and average distances of cell types to the nearest blood vessels (see Fig. \ref{fig:trivial clusters} for results in the Appendix).

\subsubsection{Multispecies witness persistence classification determines phase transitions as separate cluster}

\begin{figure}[H]
    \centering
    \captionsetup{width=1\textwidth}
        \includegraphics[width=\textwidth]{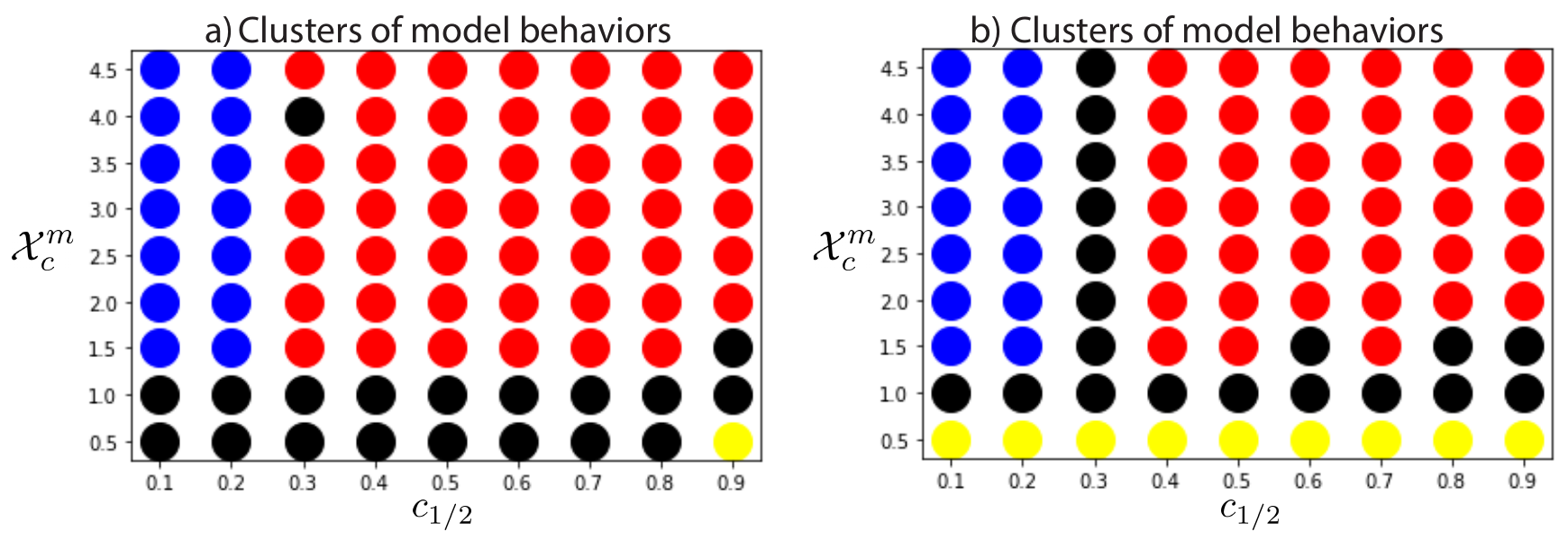}

        \caption{{\bf Classification of multispecies PH distance vectors for $k=4$ in $k$-means clustering while distinguishing between macrophage subtypes $M_1$ and $M_2$.}  a) Classification results. b) Cluster purity scores. For each parameter combination $\mathcal{X}^m_c$ and $c_{1/2}$ of the ABM, we include 20 independent simulations in our analysis. The colors red, blue, yellow, and black represent the cluster to which the majority of simulations are attributed by the $k$-means algorithm for $k=4$. 
    \label{fig:FourClusters}}
\end{figure}

Multispecies PH distance vectors further stratified the parameter space of the ABM not only into the three qualitatively different behaviors but also into the regions of phase transitions. When applying $k$-means classification for $k = 4$, the phase transitions between qualitative behaviors were identified as a separate cluster when including macrophage subtypes $M_1$ and $M_2$ in the analysis (see Fig. \ref{fig:FourClusters} b). Interestingly, when ignoring macrophage subtypes (see Fig. \ref{fig:FourClusters} a), this effect was less prominent. These results could not be obtained when using $k$-means classification for $k = 4$ on simple descriptor vectors of the data including information such as the number cells per cell type and average distances of cell types to the nearest blood vessels (see Fig. \ref{fig:trivial clusters} in the Appendix).

\subsubsection{Multispecies witness persistence classification is robust to mislabeling of cell types}

The multispecies witness PH pipeline is robust to noise introduced through relabeling. For each point cloud generated by the ABM, we relabeled up to $50\%$ of the necrotic cells, $M_1$, and $M_2$ macrophages. Relabeled cells were randomly attributed the label of one of the other two cell types. For example, a necrotic cell had a $50\%$ chance of being relabelled as a $M_1$ or $M_2$ macrophage. We focused on these three cell types because their numbers are of comparable magnitude in the ABM output, e.g., relabeling tumor cells or vessels would lead to the addition of a disproportionately high or low number of the other three cell types to the simulation output.

\begin{figure}[H]
    \centering
    \captionsetup{width=1\textwidth}
        \includegraphics[width=\textwidth]{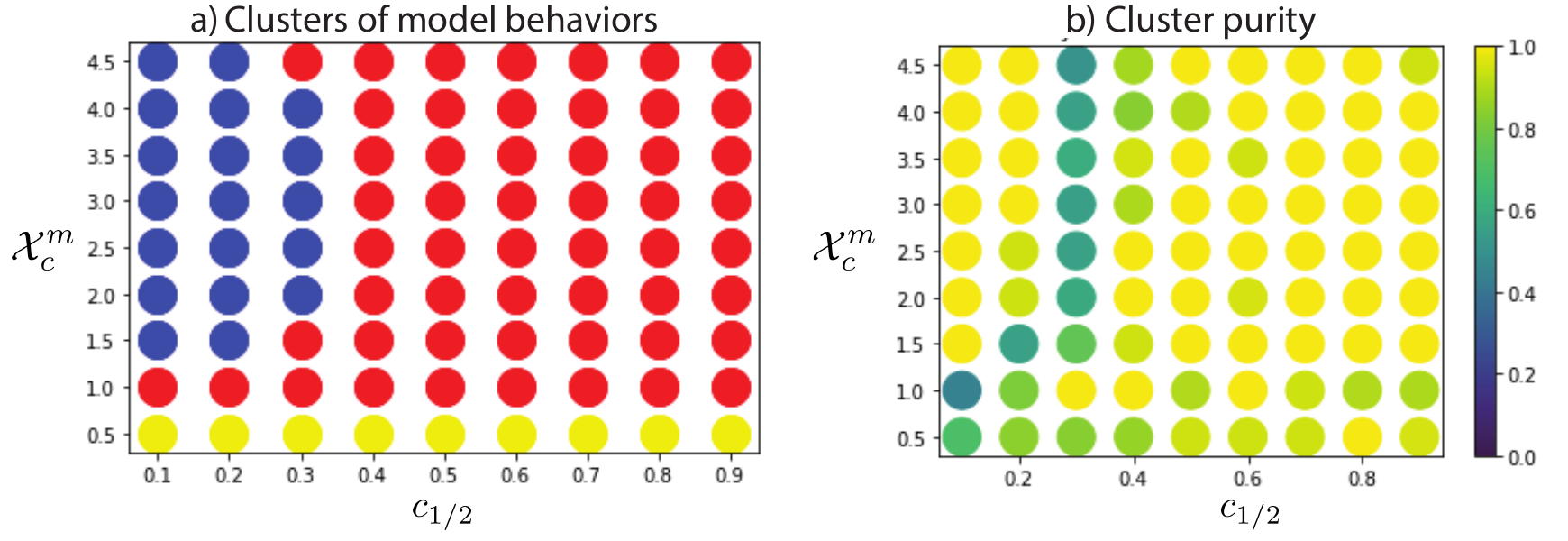}
    \caption{{\bf Classification of multispecies PH distance vectors after relabeling $50\%$ of the data.} a) Classification results. b) Cluster purity scores. For each parameter combination $\mathcal{X}^m_c$ and $c_{1/2}$ of the ABM, we include 20 independent simulations in our analysis. The colors red, blue, and yellow represent the cluster to which the majority of simulations are attributed by the $k$-means algorithm for $k=3$. The purity score is computed by taking the ratio between the number of simulations attributed to the majority clusters by 20.
    \label{fig:Noise}}
\end{figure}

\section{Discussion}
\label{sec:discussion}

With the advancement of data collection techniques, there is a growing need for analysis tools that extract relational information from spatial multispecies data. We presented two novel topological approaches to study structural relations: Dowker PH and multispecies witness PH. Dowker PH produces interpretable persistence images, but its application is limited to pairwise relations. Multispecies witness PH, on the other hand, produces features that are more difficult to interpret, but it captures relations among three or more species. We tested the utility of relational topological features in understanding macrophage and tumor behavior in point cloud simulations of the tumor microenvironment. Our results show that topological relations provide biological insight beyond that contributed by non-relational topological features and non-topological features. Furthermore, our study demonstrates that Dowker PH and multispecies witness PH effectively encode topological relations. 

This study contributes novel tools for capturing topological spatial relations that are missed in standard methods. A comparison of topological quantifications of relations to various spatial statistics \cite{WilsonMultiplex}, including the recently introduced weighted pair-correlation function \cite{Bull2023}, is postponed for future research. We believe that the topological methods, when combined with the computation of cycle representatives, may provide extra insight by identifying the local regions at which relational topological features occur in a point cloud.

Other viable topological methods include multiparameter persistence \cite{vipond2021multiparameter} and the chromatic alpha complex \cite{chromaticAlpha}. Multiparameter persistence creates multifiltrations of a simplicial complex using properties such as distances and density, and one could potentially use within-species distance and cross-species distance to create such multifiltrations. The chromatic Alpha complex creates a filtration on a multispecies version of the Delaunay triangulation. While both are viable and interesting approaches for studying multispecies data, there currently are many practical limitations to their application, such as computability and interpretability. In contrast, our approaches rely on standard one-parameter persistence, allowing efficient computation and interpretation of relational topological features.


One of the limitations of the current work is that Dowker PH can be sensitive to outliers. For example, if $P_1$ and $P_2$ are point clouds that are excluded from one another, a single outlier point of $P_1$ that lives in the neighborhood of $P_2$ will create a shared feature that is encoded by the dimension-0 Dowker persistence diagram. An enhancement of Dowker PH for robustness against outliers, possibly through subsampling \cite{StolzOutlierSubsampling, ChazalSubsampling} and multiparameter persistence, is postponed for future work.


A further study could investigate the impact of choices in the construction of multispecies witness PH. While the multispecies witness PH was based on the lazy witness complex, one could extend the construction to ``non-lazy" witness complexes. When applying the multispecies witness PH to simulated tumor microenvironments, we chose the blood vessels as landmarks. Investigation into the influence of the landmark cell type, along with the possibility of using randomly selected points in the domain as landmarks, are subjects of future work. 

Recent developments in imaging techniques \cite{CODEX_nolan, imaging_mass_cytomery} and cell identification techniques \cite{image_phenotyping_DL, celltype_morphological} produce multispecies immunohistochemistry images with detailed information about the locations of various constituents 
of a tissue microenvironment. In a tumor tissue, these constituents may include tumor cells, T-cells, B-cells, stroma, blood vessels, and more. Relational PH can potentially be applied to such multiplex images to automatically extract interpretable quantifications of relations among tumor constituents. Furthermore, the relational topological features can more broadly be applied to many other data sets which carry information on spatial locations of multiple systems. 
In the future, we envisage the integration of our methods with machine learning tools such as graph neural networks, deep learning, and random forests to achieve increased performance on such relational data and achieve novel insights.


\subsection*{Data \& code availability}
The data is available in the accompanying materials of Bull \& Byrne 2023 \cite{Bull2023}. All code is available at \url{https://github.com/irishryoon/multiplex_relations}. The Dowker PH was 
 computed in {\sc Julia} using \url{https://github.com/irishryoon/Dowker_persistence}. We implemented the multispecies witness PH in {\sc Python} using the {\sc gudhi} library \cite{maria2014} to compute persistence diagrams, as well as Bottleneck and Wasserstein distances.


\section*{Acknowledgments}
BJS, HAH, HMB, and IHRY are members of the Centre for Topological Data Analysis and this research was funded in whole or in part by EPSRC EP/R018472/1. For the  purpose of Open Access, the authors have applied a CC BY public copyright licence to any Author Accepted Manuscript (AAM) version arising from this submission. BJS is further supported by the L’Oréal-UNESCO UK and Ireland For Women in Science Rising Talent Programme. HAH gratefully acknowledges funding from EPSRC EP/K041096/1, EP/R005125/1 and EP/T001968/1, the Royal Society RGF$\backslash$EA$\backslash$201074 and UF150238, Leverhulme Trust and Emerson Collective. JAB was supported by Cancer Research UK grant number CTRQQR-2021/100002, through the Cancer Research UK Oxford Centre. IHRY gratefully acknowledges funding through the Mark Foundation for Cancer Research.

\printbibliography
\newpage
\appendix
\section{Appendix}
\subsection{Agent-based model}
\label{appendix:ABM}

We consider tumor microenvironments generated by an ABM \cite{ABM_Bonabeau}, which simulates the behavior of a system by the decisions and interactions of the agents. Our model, described in \cite{Bull2023}, simulates a growing tumor. The model explores how interactions between macrophages and the tumor microenvironment can generate interplay between varying macrophage phenotypes and the migration of tumor cells towards surrounding vasculature, a trait associated with tumor metastasis.

The model is a 2D, off-lattice, hybrid, force-based model, containing four different cell types (tumor cells, stromal cells, macrophages, and necrotic cells). Cell movement is determined by force-based interactions with neighbouring cells, together with interactions with five different chemical species described by partial differential equations (oxygen, CSF-1, TGF-$\beta$, CXCL12, and EGF). Blood vessels are represented as fixed points that are interpreted as cross-sections of vessels rising through the simulation plane. 

A key part of the model is the phenotype label $\Omega$ associated with each macrophage. The macrophage phenotype varies continuously between 0 and 1. Macrophages with $\Omega \approx 0$ are anti-tumor `$M_1$' macrophages which kill tumor cells on contact. On the other hand, `$M_2$' macrophages with $\Omega \approx 1$ are pro-tumor and produce a chemokine, EGF, which increases tumor cell migration.

Availability of oxygen mediates the cell cycle of tumor and stromal cells, with lower oxygen availability causing reduced proliferation and, with sustained lack of oxygen, death. Dead cells are labelled as necrotic, and occupy space for a period of time.
CSF-1 and CXCL12 are key chemokines for macrophages, and macrophages are attracted via chemotaxis towards increasing gradients of these chemicals. Crucially, $M_1$ macrophages are more strongly attracted towards CSF-1 (produced by tumor cells) while $M_2$ macrophages are more strongly attracted towards CXCL12 generated by perivascular fibroblasts (assumed to be co-located with blood vessels, and therefore not explicitly included as agents in the model). Macrophages enter the simulation through the vasculature with a phenotype $\Omega=0$, and are attracted towards the tumor via the CSF-1 gradient. On reaching the tumor, they are exposed to TGF-$\beta$ generated by tumor cells. Prolonged exposure to TGF-$\beta$ causes macrophage phenotype to irreversibly increase, until it reaches a maximum of $\Omega=1$. This reduces macrophage killing of tumor cells, and ultimately sensitizes them to the CXCL12 gradient produced from blood vessels, causing migration of $M_2$ macrophages back towards the vasculature. Since $M_2$ macrophages produce EGF, tumor cells can follow this gradient and may ultimately reach the vasculature (a trait associated with increased likelihood of tumor metastasis, which requires tumor cells to enter vasculature to migrate to other parts of the body).

We consider a parameter sweep in which two key parameters related to CSF-1 are varied: $\chi^m_c$, the chemotactic sensitivity of macrophages to gradients of CSF-1, and $c_{1/2}$, the concentration of CSF-1 at which macrophage extravasation is half-maximal. All other parameters are held at constant values described in \cite{Bull2023}. In Fig. \ref{fig:joshlabelsII} we show subjective classification of different qualitative behaviors of the model resulting from different parameter regimes. These qualitative behaviors manifest in different spatial distributions of the different cell types. In particular, Bull and Byrne (2023)\cite{Bull2023} relate these to the three E's of cancer immunoediting \cite{Dunn2004}: low $c_{1/2}$ leads to tumor `Elimination' as macrophages are highly recruited to the simulation and destroy the tumor. The exception to this is when $\chi_c^m$ is also low, generating `Equilibrium' behavior as macrophages are not sufficiently attracted to the tumor to destroy it. When $\chi^m_c$ and $c_{1/2}$ are both sufficiently high, $M_1$ macrophages are converted to $M_2$ macrophages faster than they can eliminate the tumor, causing tumor progression to the vasculature and thus immune `Escape'.

\newpage
\subsection{Supplementary figures and table}
\quad \newline 
\begin{figure}[h!]
    \centering
    \captionsetup{width=1\textwidth}
        \includegraphics[width=\textwidth]{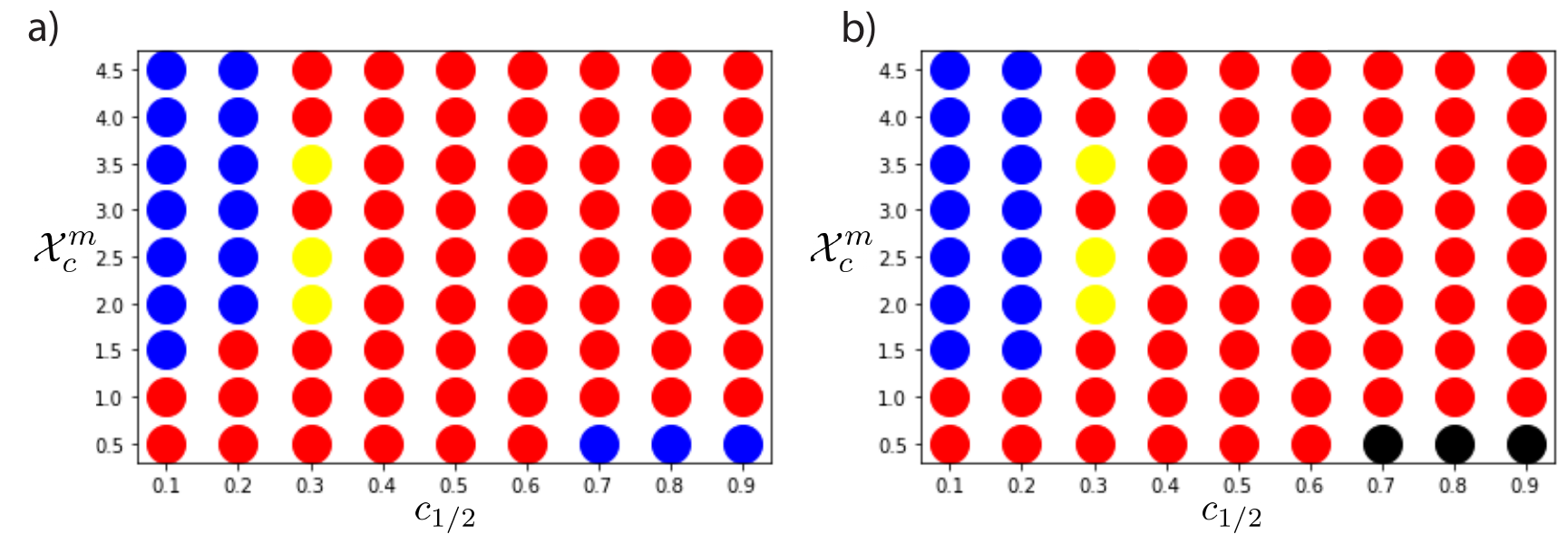}
    \caption{{\bf Classification using simple feature vectors which do not incorporate spatial information of cell types. } We apply $k$-means clustering to simple descriptor vectors for $k = 3$ (a) and $k = 4$ (b). We populate the simple description vectors with entries corresponding to the number of tumor cells, the number of macrophages, the number of necrotic cells, the average distance of tumor cells to the nearest blood vessel, the average distance of necrotic cells to the nearest blood vessel, and the average distance of macrophages to the nearest blood vessel. 
    \label{fig:trivial clusters}}
\end{figure}

\newpage

\begin{table}[ht!]
\caption{\textbf{Different versions of the multispecies witness PH and distance vectors.} We populate the multispecies distance vectors with pairwise Bottleneck $d_B$ and 1-Wasserstein distances $d_W$ for persistence diagrams $\pd_i$ in dimensions $i = 0,1$.}\label{tab:versions}
\centering
\begin{tabular}{|p{0.8cm}|p{3cm}|p{6.8cm}|}
\hline
Version & Witness filtrations considered & Distance vector entries\\
\hline
1 & tumor cells, necrotic cells, and macrophages & $ \begin{pmatrix} d_B( \pd_i(W^\bullet_{V, {T}}), \pd_i(W^\bullet_{V, {N}}))\\ d_B( \pd_i(W^\bullet_{V, {T}}), \pd_i(W^\bullet_{V, {M}})) \\ d_B( \pd_i(W^\bullet_{V, {N}}), \pd_i(W^\bullet_{V, {M}})) \\  d_W( \pd_i(W^\bullet_{V, {T}}), \pd_i(W^\bullet_{V, {N}})) \\ d_W( \pd_i(W^\bullet_{V, {T}}), \pd_i(W^\bullet_{V, {M}})) \\ d_W( \pd_i(W^\bullet_{V, {N}}), \pd_i(W^\bullet_{V, {M}})) \end{pmatrix} $
 \\
\hline
2 & tumor cells, necrotic cells, $M_1$ macrophages, and $M_2$ macrophages & 
$\begin{pmatrix}
d_B( \pd_i(W^\bullet_{V, {T}}), \pd_i(W^\bullet_{V, {N}})) \\ d_B( \pd_i(W^\bullet_{V, {T}}), \pd_i(W^\bullet_{V, {M1}})) \\ d_B( \pd_i(W^\bullet_{V, {N}}), 
\pd_i(W^\bullet_{V, {M2}})) \\ d_B( \pd_i(W^\bullet_{V, {N}}), 
\pd_i(W^\bullet_{V, {M1}}))  \\
d_B( \pd_i(W^\bullet_{V, {N}}), 
\pd_i(W^\bullet_{V, {M2}}))  \\
d_B( \pd_i(W^\bullet_{V, {M1}}), 
\pd_i(W^\bullet_{V, {M2}})) 
\\  d_W( \pd_i(W^\bullet_{V, {T}}), \pd_i(W^\bullet_{V, {N}})) \\ d_W( \pd_i(W^\bullet_{V, {T}}), \pd_i(W^\bullet_{V, {M1}})) \\ d_W( \pd_i(W^\bullet_{V, {N}}), 
\pd_i(W^\bullet_{V, {M2}})) \\ d_W( \pd_i(W^\bullet_{V, {N}}), 
\pd_i(W^\bullet_{V, {M1}}))  \\
d_W( \pd_i(W^\bullet_{V, {N}}), 
\pd_i(W^\bullet_{V, {M2}}))  \\
d_W( \pd_i(W^\bullet_{V, {M1}}), 
\pd_i(W^\bullet_{V, {M2}})) 
 \end{pmatrix} $
\\
\hline
\end{tabular}
\end{table}

\newpage

\begin{figure}[h]
    \centering
    \captionsetup{width=1\textwidth}
        \includegraphics[width=.5\textwidth]{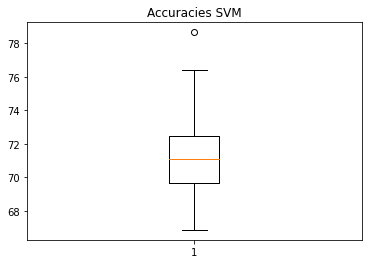}
    \caption{{\bf SVM accuracy on multispecies witness PH distance vectors.} We perform SVM analysis to predict the dominant macrophage subtype as described in Section \ref{sec:results_Dowker} on our witness feature vectors which we create as described in Section \ref{sec:WitnessFeatures} while distinguishing between $M_1$ and $M_2$ macrophages in the analysis. We present the accuracies of the SVM for 100 randomized subsets of unseen data.}
    \label{fig:SVMWitness}
\end{figure}

\newpage

\begin{figure}[H]
    \centering
    \captionsetup{width=1\textwidth}
    \includegraphics[width=\textwidth]{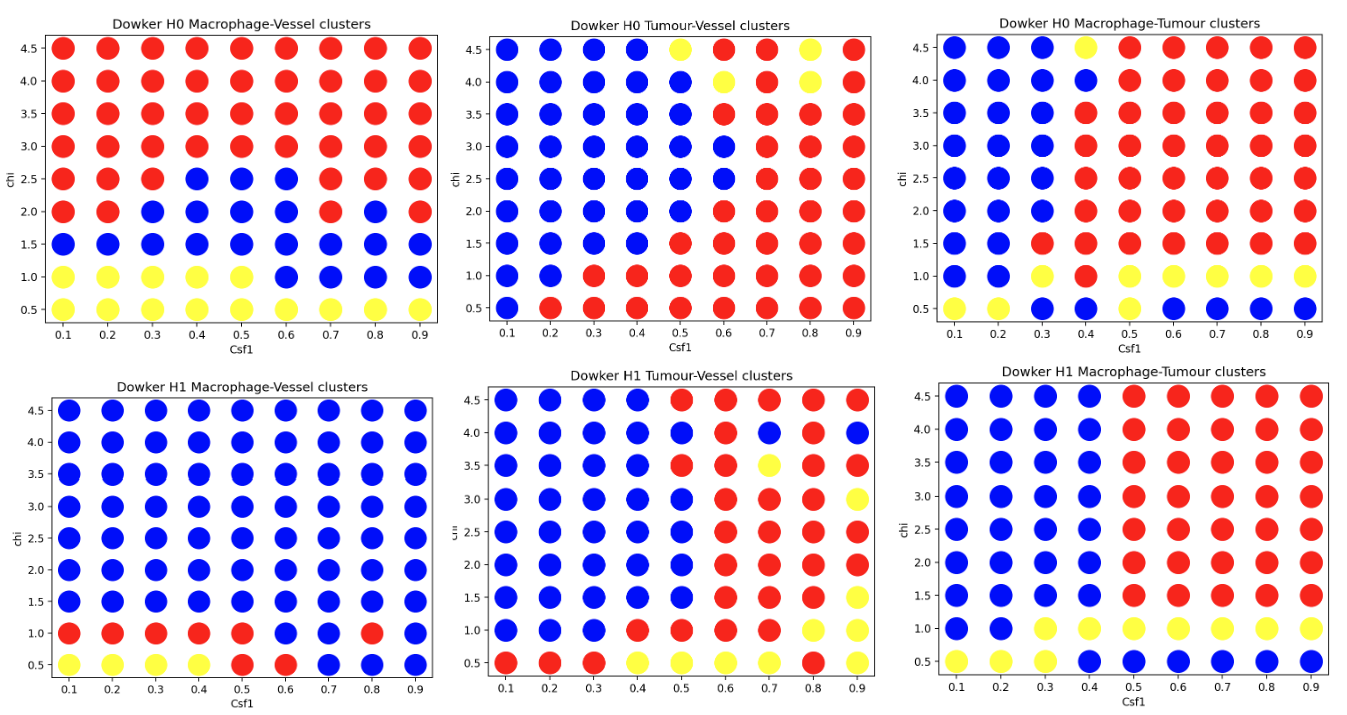}
    \caption{{\bf Classification of qualitative behavior using Dowker persistence images.} We perform clustering to infer qualitative behavior regimes as described in Section \ref{sec:WitnessFeatures} using Dowker persistence images which we create as described in Section \ref{sec:results_Dowker}. We present clustering results of each combination of tumor cells, macrophages (without knowledge of phenotype), and blood vessels in dimensions $0$ and $1$. }
    \label{fig:dowkerclustering}
\end{figure}

\end{document}